\documentclass[11pt,a4paper]{article}
\usepackage[margin=27mm]{geometry}
\usepackage[T1]{fontenc}
\usepackage{lmodern}
\usepackage{amsmath,amssymb,amsthm,mathtools}
\usepackage[expansion=false]{microtype}
\usepackage[colorlinks=true,linkcolor=blue,citecolor=blue,urlcolor=blue]{hyperref}
\usepackage{mathrsfs}
\numberwithin{equation}{section}
\newtheorem{theorem}{Theorem}[section]
\newtheorem{lemma}[theorem]{Lemma}
\newtheorem{proposition}[theorem]{Proposition}
\newtheorem{corollary}[theorem]{Corollary}
\theoremstyle{definition}
\newtheorem{definition}[theorem]{Definition}

\newcommand{\R}{\mathbb R}
\newcommand{\Z}{\mathbb Z}

\DeclareMathOperator{\SL}{SL}
\DeclareMathOperator{\NE}{NE}
\DeclareMathOperator{\NUE}{NUE}
\DeclareMathOperator{\area}{area}
\DeclareMathOperator{\covol}{covol}
\DeclareMathOperator{\diag}{diag}
\DeclareMathOperator{\inj}{inj}
\DeclareMathOperator{\dist}{dist}
\DeclareMathOperator{\vol}{vol}

\allowdisplaybreaks[1]
\title{Hausdorff dimension of nonergodic directions in higher genus via effective equidistribution}
\author{Yuming Wei}
\date{}
\begin{document}
\maketitle

\begin{abstract}
For each fixed $g\ge3$, we prove that non-uniquely ergodic directions have Hausdorff dimension $1/2$ on almost every area-one surface in $\mathcal H^{\mathrm{odd}}(g-1,g-1)$ when $g$ is odd, and in $\mathcal H^{\mathrm{nonhyp}}(g-1,g-1)$ when $g$ is even. We construct these directions through cyclic splittings into $g$ slit tori. Effective equidistribution of Lindenstrauss--Mohammadi--Yang supplies sufficiently many children at each normal input. We estimate cusp and periodic-orbit failures with the original branching weights and prove that their total mass loss is summable. The surviving path measure yields the dimension lower bound, also for directions nonergodic for flat area.
\end{abstract}

\section{Introduction}\label{sec:intro}
A translation surface of genus \(g\ge2\) is a pair \((X,\omega)\), where \(X\) is a compact Riemann surface of genus \(g\) and \(\omega\) is a nonzero holomorphic one-form. Let \(p_1,\ldots,p_s\) be the zeros of \(\omega\), with respective orders \(k_1,\ldots,k_s\). The flat metric induced by \(\omega\) has cone angle \(2\pi(k_i+1)\) at \(p_i\), and the zero orders satisfy
$\sum_{i=1}^{s}k_i=2g-2$.
The stratum \(\mathcal H(k_1,\ldots,k_s)\) parametrizes isomorphism classes of translation surfaces with these zero orders. We work with surfaces of unit flat area,
\[
\operatorname{Area}(X,\omega)
=\frac{i}{2}\int_X\omega\wedge\overline{\omega}
=1,
\]
and denote the corresponding locus by \(\mathcal H_1(k_1,\ldots,k_s)\).
In genus two, the only strata are \(\mathcal H(2)\) and \(\mathcal H(1,1)\), both of which are connected. In higher genus, strata need not be connected. In particular, the classification of Kontsevich and Zorich \cite[Theorems~1 and~2]{KontsevichZorich} implies that \(\mathcal H(g-1,g-1)\) is disconnected for every \(g\ge3\). For each such \(g\), we focus on the connected component
\begin{equation}\label{eq:components}
\mathcal C_g=
\begin{cases}
\mathcal H_1^{\mathrm{odd}}(g-1,g-1),&g\text{ odd},\\
\mathcal H_1^{\mathrm{nonhyp}}(g-1,g-1),&g\text{ even}.
\end{cases}
\end{equation}
Here \(\mathcal H_1^{\mathrm{odd}}(g-1,g-1)\) denotes the odd-spin component distinct from the hyperelliptic component, while \(\mathcal H_1^{\mathrm{nonhyp}}(g-1,g-1)\) denotes the nonhyperelliptic component in even genus.

For a direction
$\theta\in \mathbb{RP}^1$, let $\phi_\theta^s$ denote its straight-line flow,
defined for almost every point and every fixed time $s$. We say $\theta$ is an ergodic direction if the flow $\phi_\theta^s$ is ergodic with the normalized flat area measure $\mu_X$. If $\mu_X$ is the only invariant ergodic measure, we say $\theta$ is a uniquely ergodic direction.
Write
$\NE(X,\omega)$ for the directions in which this flow is not ergodic
with respect to flat area. Write $\NUE(X,\omega)$ for the directions in which
the flow admits more than one invariant Borel probability measure
on its nonsingular full-orbit domain.  We work with
the measure on area-one connected components 
and show the following result:
\begin{theorem}\label{thm:main}
For every fixed $g\ge3$, for Masur--Veech almost every $(X,\omega)\in\mathcal C_g$,
\begin{equation}\label{eq:main-dimension}
\dim_H\NE(X,\omega)=\dim_H\NUE(X,\omega)=\frac12.
\end{equation}
\end{theorem}

Masur's upper bound gives $\dim_H\NUE(X,\omega)\le1/2$ \cite{Masur1992}. 
When $g=2$, $\mathcal H_2$ has two strata $\mathcal H(2)$ and $\mathcal H(1,1)$ where both are connected.
Athreya--Chaika \cite{AthreyaChaika} first showed that for almost every translation surfaces in $\mathcal H(2)$, the set of non-uniquely ergodic directions has 1/2 Hausdorff dimension. Chaika and Masur \cite{ChaikaMasur} 
subsequently proved the corresponding equality for almost every surface in the hyperelliptic components
\begin{equation*}
\mathcal H_1^{\mathrm{hyp}}(2g-2)
\quad\text{and}\quad
\mathcal H_1^{\mathrm{hyp}}(g-1,g-1),
\qquad g\geq2.
\end{equation*}
In particular, their theorem includes the stratum $\mathcal H_1(1,1)$. The proofs of Athreya–Chaika and Chaika–Masur use Rauzy induction for interval exchange transformations.

Theorem 1.1 establishes the same sharp value in the nonhyperelliptic components $\mathcal C_g$ defined in \eqref{eq:components}: the odd-spin component of $\mathcal H_1(g-1,g-1)$ for odd $g\geq3$, and its nonhyperelliptic component for even $g\geq4$. Together with Chaika–Masur’s theorem, this gives
$$\dim_H\operatorname{NUE}(X,\omega)=\frac12$$
for Masur–Veech almost every surface in every connected component of ($\mathcal H_1(g-1,g-1)$) when $g=3$ or $g\geq4$ is even; see Corollary~\eqref{cor:components}. For odd ($g\geq5$), the present argument treats the odd-spin component and does not address the even-spin component. Moreover, our construction directly produces directions nonergodic with respect to flat area, giving the equality for both $\NE(X,\omega)$ and $\NUE(X,\omega)$ on almost every surface in $\mathcal C_g$.

Slit constructions also play an important role in earlier dimension results. For the classical one-parameter family of genus-two surfaces obtained as branched double covers of a flat torus, Cheung\cite{Cheung} proved the sharp dimension bound under a Diophantine condition on the slit parameter. Cheung–Hubert–Masur\cite{CheungHubertMasur} subsequently established a $0$ or $1/2$ dichotomy characterized by a continued-fraction condition. The construction considered here allows the torus periods to vary independently in an open family of higher-genus surfaces and combines simultaneous lattice counting with estimates for the mass lost over successive generations. More recently, Wei–Yang\cite{WeiYang} developed a quantitative splitting construction in the eigenform loci $\mathcal E_D\subset\mathcal H(1,1)$ for nonsquare discriminants $D$. They proved that the set of non-uniquely ergodic directions has Hausdorff dimension $1/2$ on every surface in these loci, except for the Teichmüller curve generated by the regular decagon in $\mathcal E_5$. Their argument combines primitive lattice-point estimates with a selection among three twist parameters to obtain sufficiently many children preserving a normality condition on the splittings.

The present paper uses the cyclic slit construction of Eskin–Masur–Zorich \cite[\S9.3]{EskinMasurZorich}, specialized to $g$ flat tori. Section \S\ref{sec:cyclic} shows that these splittings form an open family in $\mathcal C_g$. We iterate twists that preserve the cyclic decomposition, choosing one loop shorter than the other $g-1$ loops so that the new slits remain embedded and the changes in the torus areas are summable. The admissible choices are encoded by primitive-vector windows in a product of $g$ lattice spaces. Effective equidistribution of Lindenstrauss–Mohammadi–Yang \cite{LindenstraussMohammadiYang} then supplies quantitative branching estimates at inputs satisfying the normality conditions of Section \S\ref{sec:normal}. A related approach using effective equidistribution for two-torus splittings is developed in \cite{Wei}.

The main difficulty is to control the cumulative effect of inputs at which these branching estimates are unavailable. The distributions of later splittings depend on the preceding choices, so estimates for exceptional sets with respect to ambient volume alone do not control the mass lost during the iteration. We retain the original equal branching weights and average the incoming mass over a positive-volume family of initial splittings. \S\ref{sec:fibers} and \S\ref{sec:survival} balance counts of integral histories against volumes in period coordinates to estimate the mass first discarded at each depth because of cusp or periodic obstructions. These losses are summable, and their total tends to zero as the initial slit height tends to infinity. Consequently, a positive-measure family of initial splittings supports a path measure of positive surviving mass. Summable area exchange yields nonergodicity along every surviving path, while the original cylinder weights give the Hausdorff dimension lower bound. Finally, $\mathrm{SL}_2(\mathbb R)$-invariance and ergodicity of Masur–Veech measure promote the conclusion to almost every surface in ($\mathcal C_g$).

We record the quantitative statement used in this argument. Fix $r=1+t>1$ and $\varepsilon>0$. At slit height $W$, put
\begin{equation}\label{eq:scales}
\gamma=r(1+\varepsilon),\qquad V=W^r,\qquad U=W^\gamma,
\qquad \delta=\frac{d}{\log W},\qquad \tau=\frac VU,\qquad Q=\frac UW,
\end{equation}
where $d>0$ is fixed and small. Constants may depend on the fixed genus, area margin and parameters, but not on height or depth; this convention also applies to $\ll$ and $\asymp$. Throughout, $x\times y=\det(x,y)$, $J(x,y)=(-y,x)$, and $\lambda_1(\Lambda)$ is the length of a shortest nonzero vector of a unimodular lattice. Angles are the smaller positive angles between vectors in the common narrow cone of the construction.

\begin{theorem}\label{thm:quantitative}
For every fixed $g\ge3$ and every sufficiently small $t>0$, there are arbitrarily small $\varepsilon>0$ for which the following holds. There is a positive-volume family of tall cyclic roots and a Borel normality test, given in Definition~\ref{def:normal}, such that each normal input $S$ of height $W$ in the resulting trees has at least
\begin{equation}\label{eq:branching}
b(S)\ge c_g\frac{W^t}{(\log W)^g}
\end{equation}
distinct children. Their slit lengths are comparable to $W^\gamma$, their directions are separated by $c_gW^{-2\gamma}$, and every component area remains uniformly positive along every history.

Give each child the fraction $1/b(S)$ of its parent's incoming mass. Discard a child when its input is not normal, keeping all other weights unchanged. If $L_n$ is the root-averaged mass first discarded at depth $n$, then, for some $c_*>0$,
\begin{equation}\label{eq:mass-loss}
L_n\le \exp\{-c_*\gamma^n\log W_0\}\quad(n\ge1),
\qquad L_0+\sum_{n\ge1}L_n\longrightarrow0\quad(W_0\longrightarrow\infty).
\end{equation}
Consequently, a positive-measure set of roots has positive surviving mass.
\end{theorem}

\textbf{AI assistance.} The idea of using effective equidistribution to obtain quantitative branching estimates for splittings of translation surfaces in $\mathcal H(1,1)$ arose in April 2026 during joint work on establishing Hausdorff dimension $1/2$ for non-uniquely ergodic directions in eigenform loci with Pengyu Yang. The author subsequently extended this approach to higher genus. ChatGPT 6.1 sol was used to assist in checking the calculations and estimates for bad inputs in the period coordinates of the initial splittings, and to proofread and improve the exposition. The author takes responsibility for the mathematical content and the final text.

\section{Connected components and cyclic splittings}\label{sec:cyclic}
Let $\Omega(X)$ denote the complex vector space of holomorphic
$1$-forms on a compact Riemann surface $X$ of genus $g\ge2$.
The surface $X$ is \emph{hyperelliptic} if it admits a
holomorphic involution $\iota$ satisfying
$\iota^*\eta=-\eta$ for every $\eta\in\Omega(X)$;
equivalently, $X/\langle\iota\rangle$ is the Riemann sphere.
Otherwise, $X$ is \emph{nonhyperelliptic}.
In $\mathcal H(g-1,g-1)$, the hyperelliptic component
$\mathcal H^{\mathrm{hyp}}(g-1,g-1)$ consists of holomorphic
$1$-forms
on hyperelliptic curves whose hyperelliptic involution
interchanges the two zeros
\cite[\S2.1]{KontsevichZorich}.

Suppose that the zeros $p_1,\ldots,p_s$ of $\omega$ have even
orders $2m_1,\ldots,2m_s$.
Let $V_\omega\subseteq\Omega(X)$ be the subspace of holomorphic
$1$-forms vanishing to order at least $m_i$ at each $p_i$.
The \emph{spin parity} of $(X,\omega)$ is
$\phi(X,\omega)=\dim_{\mathbb C}V_\omega\bmod2$.
The surface has \emph{even spin} or \emph{odd spin} according
as this value is $0$ or $1$.
This parity is constant on each connected component of the
stratum \cite[\S2.2]{KontsevichZorich}.

We first recall some results that we will use later.
\begin{theorem}[Period coordinates,{\cite[\S1.1]{KontsevichZorich}}]
\label{thm:period-coordinates}
Let $g\ge2$. After passing to a marked local chart, integration of the differential
over an integral basis of $H_1(X,\Sigma;\Z)$ gives local complex
coordinates on $\mathcal H(k_1,\ldots,k_s)$. Their number is
$2g+s-1$. Lebesgue measure in these coordinates induces the
Masur--Veech measure class, and radial disintegration defines its
area-one measure class.
\end{theorem}

\begin{theorem}[Kontsevich--Zorich, {\cite[Theorems 1 and 2]{KontsevichZorich}}]
\label{thm:KZ-specialized}
The components of $\mathcal H(g-1,g-1)$ are as follows.
For $g=2$ the stratum is connected and hyperelliptic; for $g=3$
its components are hyperelliptic and odd-spin. For even $g\ge4$
they are hyperelliptic and nonhyperelliptic. For odd $g\ge5$
they are hyperelliptic, even-spin and odd-spin, with the last two
nonhyperelliptic. Spin parity is constant on each component where
it is defined.
\end{theorem}

\begin{theorem}[Eskin--Masur--Zorich,{\cite{EskinMasurZorich}}]
\label{thm:EMZ-slits}
If the resulting
surface and its closed component surfaces have only even zeroes,
their spin parities satisfy
\begin{equation}\label{eq:spin-additivity}
 \varphi(X,\omega)=\sum_i\varphi(X_i,\omega_i)\pmod2.
\end{equation}
This is \cite[Lemma 10.1]{EskinMasurZorich}. Moreover, for almost
every surface in $\mathcal H^{\mathrm{hyp}}(g-1,g-1)$, saddle
connections between the distinct zeroes have multiplicity at most
two \cite[Corollary 10.4]{EskinMasurZorich}.
\end{theorem}
Let \(T_i=\mathbb R^2/L_i\), \(i=1,\ldots,g\), be flat tori carrying embedded slits \(\ell_i\) with a common nonzero holonomy vector \(w\). Cut each \(T_i\) along \(\ell_i\), and identify the positive bank of \(\ell_i\) with the negative bank of \(\ell_{i+1}\) by translation, matching the initial and terminal endpoints, with indices taken modulo \(g\). This produces a translation surface \((X,\omega)\in\mathcal H(g-1,g-1)\). We call the resulting decomposition a \emph{cyclic splitting} of \((X,\omega)\).
Let \(\sigma_i\) denote the saddle connection arising from the identification between \(T_i\) and \(T_{i+1}\). We refer to \(\sigma_1,\ldots,\sigma_g\) as the \emph{slits of the splitting}. They have pairwise disjoint interiors and common holonomy \(w\) when oriented from the common initial endpoint to the common terminal endpoint. Cutting along all these slits recovers the \(g\) slit tori.

Write \(A_i=\operatorname{covol}(L_i)\) and normalize the total area so that \(\sum_{i=1}^g A_i=1\). Fix \(0<a<1/(2g)\). All local estimates below are made under the bounds \(A_i\ge a\) for every \(i\).
This construction is the special case of Theorem~\ref{thm:EMZ-slits} in which every closed component surface is a torus. To apply its short-slit description, keep the lattices and the initial endpoints fixed and shorten each \(\ell_i\) along its original direction, replacing its holonomy \(w\) by \(\rho w\), where \(0<\rho\le1\). The shortened slits remain embedded, and the resulting surfaces form a continuous family in \(\mathcal H(g-1,g-1)\), with the original surface corresponding to \(\rho=1\). Hence the entire family lies in the same connected component. Consequently, the connected-component conclusions obtained for sufficiently small \(\rho>0\) also hold for the original surface; when \(g\) is odd, the same applies to the spin parity.

\begin{lemma}\label{lem:component}
For $g\ge3$, every surface obtained by the cyclic construction belongs to $\mathcal C_g$. The torus periods together with $w$ give local relative period coordinates; after area normalization, the image is locally open in $\mathcal C_g$.
\end{lemma}

\begin{proof}
The initial endpoints become one zero and the terminal endpoints another, each with cone angle $2\pi g$. Both zero orders are therefore $g-1$. The boundary of each cut torus is the difference of two adjacent slits, so the $g$ oriented slits represent the same relative homology class.

Choose a symplectic pair of cycles on each torus, away from the slit. The resulting $2g$ cycles have the standard unimodular intersection matrix and hence form an integral symplectic basis of absolute homology. Adjoining one slit gives a relative homology basis. Their periods are exactly the torus periods and $w$, and small changes preserve the embedded slits and the gluing. The period map is therefore a local coordinate chart. Restriction to unit area remains locally open, since radial dilation is transverse to the area-one locus.

The $g$ distinct homologous slits persist throughout this open family. If it met the hyperelliptic component, that component would contain a positive-measure open set with at least three homologous saddle connections joining its two zeros. This contradicts the almost-everywhere multiplicity bound of two in Theorem~\ref{thm:EMZ-slits}.

For odd $g$, \eqref{eq:spin-additivity} gives
$\phi(X,\omega)=\sum_{i=1}^g\phi(T_i)=g=1\pmod2$,
because each torus has odd spin. The surface is therefore in the odd-spin component. For even $g$, exclusion of the hyperelliptic component leaves the nonhyperelliptic component. The by Theorem~\ref{thm:KZ-specialized} gives \eqref{eq:components}.
\end{proof}

For $g=2$, the same gluing lies in the connected hyperelliptic stratum $\mathcal H(1,1)$; the corresponding dimension result is known \cite{ChaikaMasur,Wei}. The quantitative estimates below assume $g\ge3$.

\section{Dehn twists on the cyclic decomposition}\label{sec:move}
The next Lemma introduces a Dehn twist operation on the translation surface.
\begin{lemma}\label{lem:move}
Let $v_i\in L_i$ be primitive and satisfy $0<D_i:=w\times v_i<A_i$ for every $i$. Put
$v_\Sigma=\sum_{i=1}^g v_i$, and $w^+=w+v_\Sigma$.
If $D_i^+:=w^+\times v_i>0$ for every $i$, there is a new cyclic splitting of the same translation surface with common slit $w^+$. For any positively oriented basis $(v_i,u_i)$ of $L_i$, its lattices and areas are
\begin{equation}\label{eq:move}
L_i^+=\Z v_i+\Z\left(u_i-\sum_{j\ne i}v_j\right),
\qquad A_i^+=A_i-v_i\times v_\Sigma.
\end{equation}
The old and new torus regions satisfy
\begin{equation}\label{eq:area-exchange}
\sum_{i=1}^g\area(T_i\triangle T_i^+)\le2\sum_{i=1}^gD_i.
\end{equation}
\end{lemma}

\begin{proof}
The $v_i$-parallel leaves missing the slit form a cylinder $C_i$ of area $A_i-D_i$. After cutting along the slit, the complementary region develops to a parallelogram with sides $v_i,w$. Glue these parallelograms along their slit sides, retaining distinct boundary vertices until the cylinders are reattached. The result is a flat annulus whose interior embeds in the original surface. We retain the $C_i$ and change only the cuts in this annulus.

Rotate so that $w=(h,0)$ with $h=|w|$, and write $v_i=(x_i,y_i)$, where $y_i=D_i/h>0$. For the development, order the indices as encountered in the positive loop directions; this may reverse the original cyclic order and does not change any formula involving the individual loops and their sum. Set
\[
P_0=0,\qquad P_i=\sum_{j=1}^i v_j,\qquad P_g=(p,q)=v_\Sigma,
\qquad P_{i+g}=P_i+(p,q).
\]
The universal cover develops to the polygonal strip
\begin{equation*}
B=\{(x,y)\in\R^2:f(y)\le x\le f(y)+h\},
\end{equation*}
where $f$ is the piecewise affine function whose graph passes through the $P_i$ and satisfies $f(y+q)=f(y)+p$. Its successive slopes are $x_i/y_i$, and the deck translation is $(x,y)\mapsto(x+p,y+q)$.

Put $m=(h+p)/q$. For $0\le i<g$, consider the segment
$\gamma_i(t)=P_i+(mt,t)$ for $0\le t\le q$,
from $P_i$ to $P_i+w^+$. Its horizontal distance from the left boundary is
$d_i(t)=P_i^x+mt-f(P_i^y+t)$.
Away from boundary vertices,
\begin{equation*}
d_i'(t)=m-\frac{x_j}{y_j}=\frac{D_j^+}{qy_j}>0,
\qquad d_i(0)=0,\qquad d_i(q)=h.
\end{equation*}
Thus the interior of each segment lies in the interior of $B$. In the transverse coordinate $s=x-my$, consecutive starting values differ by $-D_i^+/q<0$, and the deck translation changes $s$ by $p-mq=-h$. The $g$ starting values are distinct modulo $h$. Consequently, the projected segments are embedded and pairwise disjoint. They divide the annulus into parallelograms with sides $v_i,w^+$, whose boundary edges attach to the same cylinders $C_i$. Each new cut travels one period further around the annulus, so the replacement is a Dehn twist of the splitting.

A connector through $C_i$, closed by the reversed old slit, represents the complementary period $u_i$ and has holonomy $u_i+w$. Closing it instead by the reversed new slit gives $u_i-v_\Sigma$. The unchanged cylinder core and this closed connector have intersection number one, so their classes form a homology basis of the new torus. Adding $v_i$ to the second generator proves \eqref{eq:move}, and
\[
A_i^+=(A_i-D_i)+D_i^+=A_i-v_i\times v_\Sigma>0.
\]
Both torus regions contain $C_i$, so their symmetric difference has area at most $D_i+D_i^+$. Since
\[
\sum_iD_i^+=w^+\times v_\Sigma=w\times v_\Sigma=\sum_iD_i,
\]
we obtain \eqref{eq:area-exchange}. The common slit holonomy and cyclic adjacency are preserved, giving the required splitting of the original surface.
\end{proof}

We now choose windows ensuring the positivity in Lemma~\ref{lem:move}. Rotate $w$ to $(W,0)$ and set
\begin{equation}\label{eq:moving-frame}
\begin{gathered}
G=\diag\left(\sqrt{\frac{\delta}{WU}},\sqrt{\frac{WU}{\delta}}\right),
\qquad D_Q=\diag(Q^{-1/2},Q^{1/2}),\\
\ell_b=\begin{pmatrix}1&0\\ b&1\end{pmatrix},
\qquad H_b=D_Q\ell_{-b}G.
\end{gathered}
\end{equation}
Fix $\kappa=1/(100g)$ and $\alpha_i=1+2(i-2)\kappa$ for $2\le i\le g$. Choose $K=K(g)\ge10g$ sufficiently large and then $c=c(g)>0$ sufficiently small, with $c<\kappa/100$. Recall that $\tau=\frac VU$, 
and for $b\in I:=[K-1,K+1]$, use the target boxes
\begin{equation}\label{eq:boxes}
\begin{aligned}
B_1&=\{\sqrt\delta(x,y):\tau<x<2\tau,\ 1<y<2\},\\
B_i&=\{\sqrt\delta(x,y):\alpha_i<x<\alpha_i+\kappa,\ |y|<c\},
\qquad 2\le i\le g.
\end{aligned}
\end{equation}
All constants in this section may depend on the
fixed genus and the positive lower bound for the areas.

\begin{definition}[Children of splitting]
A child of a splitting $(L_1,...,L_g,w)$ is a slit $w+\sum v_i$ satisfying the following conditions:\\
(i) $(v_1,\ldots,v_g)$ is a primitive tuple in Lemma~\ref{lem:move};\\
(ii) $H_bv_i\in B_i$ for every $i$ at some common $b\in I$.
\end{definition}

\begin{lemma}\label{lem:children}
For sufficiently large $W$, the children $w^+=w+\sum v_i$ share the following properties
\begin{equation*}
\begin{gathered}
|v_1|\asymp V,\qquad |v_i|\asymp U\quad(2\le i\le g),\qquad |w^+|\asymp U,\\
0<w\times v_i\le C\delta,\qquad c_1\delta\le w^+\times v_i\le C_1\delta,\\
\angle(w,v_i)\le\frac{C\delta}{WU},\qquad
\angle(w,w^+)\le\frac{C\delta}{WU}.
\end{gathered}
\end{equation*}
Distinct children give distinct directions that are separated by at least $c_2U^{-2}$.
\end{lemma}

\begin{proof}
Write $H_bv_i=\sqrt\delta(\xi_i,\eta_i)$. Inverting \eqref{eq:moving-frame} gives
\begin{equation}\label{eq:window-coordinates}
(v_i)_x=U\xi_i,\qquad
(v_i)_y=\frac{\delta}{W}\left(b\xi_i+\frac WU\eta_i\right),
\qquad w\times v_i=\delta\left(b\xi_i+\frac WU\eta_i\right).
\end{equation}
Since $(W/U)/\tau=W/V\to0$, the coordinate ranges in \eqref{eq:boxes} imply
$w\times v_1\asymp K\tau\delta$ and $w\times v_i\asymp K\delta$ for $i\ge2$. Determinant preservation gives, for large $W$,
\begin{equation*}
-5\delta\le v_1\times v_i\le-\tfrac12\delta\quad(i\ge2),
\qquad |v_i\times v_j|\le4c\delta\quad(i,j\ge2).
\end{equation*}
Hence $w^+\times v_1\ge(g-1)\delta/2$, while for $i\ge2$ the contribution $w\times v_i\ge(K-2)\delta$ dominates the possible negative contribution $(5+4gc)\delta$. All new cross products are therefore positive and comparable to $\delta$. Also $w\times v_i<A_i$ for large $W$, so Lemma~\ref{lem:move} applies.

The length estimates follow from the positive horizontal coordinates in \eqref{eq:window-coordinates}. The short-loop angle is at most
$\frac{C\delta}{WU}+\frac{C\delta}{UV}\le\frac{C'\delta}{WU}$,
and the other loop angles satisfy the same bound. The horizontal coordinate of $w^+$ is $W+U\sum_i\xi_i\asymp U$, and its vertical coordinate is $O(\delta/W)$, giving the remaining estimates.

Fix $(v_1,v_2)$. For $i\ge3$, write a possible partner as
\[
v_i=a_iv_2+e_i\frac{Jv_2}{|v_2|^2},
\qquad a_i\asymp1,\qquad e_i=v_2\times v_i=O(c\delta).
\]
For two possible choices $v_i,v_i'$ in $L_i$, even with different witnessing times, their determinant is $a_ie_i'-a_i'e_i=O(c\delta)$. But $|v_i\times v_i'|$ is an integer multiple of $A_i$ and hence zero for large $W$. Positivity and primitivity then imply $v_i=v_i'$. Thus $(v_1,v_2)$ determines every remaining loop.

Distinct positive primitive choices of $v_1$ have directions separated by at least $cV^{-2}$, since their nonzero determinants have magnitude at least $A_1\ge a$. Their children lie within $C\delta/(VU)$ of the respective $v_1$ directions. Because $\delta V/U\to0$, these children remain separated by $c'V^{-2}$. With $v_1$ fixed, distinct choices of $v_2$ are separated by at least $cU^{-2}$, whereas each child lies within $C\delta/U^2$ of its $v_2$ direction. Since $\delta\to0$, their children remain separated by $c'U^{-2}$. This proves both distinctness and the claimed separation.
\end{proof}

\section{A homogeneous dynamical model}\label{sec:orbit}
Put $X_1=\SL_2(\R)/\SL_2(\Z)$ and $X=X_1^g$, with product Haar probability measure $m_X$. We identify $h\SL_2(\Z)$ with $h\Z^2$. The primitive Siegel formula \cite{Siegel}, with both signs of every primitive vector included, is
\begin{equation}\label{eq:siegel}
\int_{X_1}\sum_{v\in\Lambda^{\mathrm{prim}}}f(v)\,dm_1(\Lambda)
=\frac1{\zeta(2)}\int_{\R^2}f(x)\,dx
\end{equation}
for nonnegative Borel $f$. In particular,
\begin{equation*}
m_1\{\Lambda:\lambda_1(\Lambda)<\epsilon\}\le\frac{\pi\epsilon^2}{\zeta(2)}.
\end{equation*}

For the homogeneous model, express the slit $w$ and the lattices $L_i$
in fixed period coordinates. Let $R_\theta$ denote counterclockwise
rotation through angle $\theta$, and put $W=|w|$. We use the
normalization
\begin{equation}\label{eq:Lambda_i}
\begin{gathered}
G_w=\diag\left(
\sqrt{\frac{d}{W^{1+\gamma}\log W}},
\sqrt{\frac{W^{1+\gamma}\log W}{d}}\right)R_{-\arg w},\\
\Lambda_i=A_i^{-1/2}G_wL_i,\qquad 1\le i\le g.
\end{gathered}
\end{equation}
Since $\det G_w=1$ and $\covol(L_i)=A_i$, each $\Lambda_i$ is
unimodular. All subsequent occurrences of the input lattices
$\Lambda_i$ use this normalization.

Since $U=W^\gamma$ and $\delta=d/\log W$, the diagonal factor in
$G_w$ is the matrix $G$ in \eqref{eq:moving-frame}. Set
$H_{w,b}:=D_Q\ell_{-b}G_w$. Then
$H_{w,b}v=H_b(R_{-\arg w}v)$, so the window conditions of
Section~\ref{sec:move} take the form $H_{w,b}v_i\in B_i$ for
$v_i\in L_i^{\mathrm{prim}}$ at a common $b\in I$.
Lemma~\ref{lem:children} applies to these conditions by rotating
the entire splitting to make the slit horizontal.

Choose smooth $0\le\psi_i\le1$ supported in $B_i$ and equal to one
on fixed-proportion smaller boxes. Set
\begin{equation}\label{eq:test-functions}
F_i(M)=\sum_{v\in M^{\mathrm{prim}}}\psi_i(\sqrt{A_i}\,v),
\qquad \Phi(M_1,\ldots,M_g)=\prod_{i=1}^gF_i(M_i).
\end{equation}
Two vectors in the closure of one box have determinant smaller than $a\le A_i$. Thus independent lattice vectors cannot both contribute, and parallel primitive vectors in a positive longitudinal window coincide. Each sum has at most one nonzero term, so $0\le\Phi\le1$. Applying \eqref{eq:siegel} in each factor and Tonelli's theorem gives
\begin{equation}\label{eq:test-mean}
\int_X\Phi\,dm_X
=\prod_{i=1}^g\frac{\int_{\R^2}\psi_i}{A_i\zeta(2)}
\asymp_g\delta^g\tau.
\end{equation}

Write $a_s=\diag(e^s,e^{-s})$, $u_b=\bigl(\begin{smallmatrix}1&b\\0&1\end{smallmatrix}\bigr)$ and $k=J$. A matrix without a factor index acts diagonally on $X$.

\begin{lemma}\label{lem:test}
Let $s=\tfrac12\log\delta^{-1}$, $h_\delta=a_su_{-1}$ and $\Psi=\Phi\circ k^{-1}\circ h_\delta^{-1}$. There is a compact set $\mathcal K\subset X$, independent of $W$ and of the component areas in their prescribed compact positive range, such that $\operatorname{supp}\Psi\subset\mathcal K$. For every fixed finite-order Sobolev norm $S$ on $X$, there are $P,d_1,C$, depending only on the fixed data, such that
\begin{equation}\label{eq:sobolev}
S(\Psi)\le C\delta^{-P}\tau^{-d_1}.
\end{equation}
With
\begin{equation}\label{eq:orbit-input}
T=\tfrac12\log(Q/\delta),\qquad
z=u_{K-1/Q}k(\Lambda_1,\ldots,\Lambda_g),
\end{equation}
we have
\begin{equation}\label{eq:orbit-identity}
\Psi(a_Tu_bz)=\Phi\bigl(D_Q\ell_{-(K+b)}(\Lambda_1,\ldots,\Lambda_g)\bigr),
\qquad -1\le b\le1.
\end{equation}
\end{lemma}

\begin{proof}
For a  vector $x=(x_1,x_2)\in B_i$,
\[
h_\delta kx=\left(-\frac{x_1+x_2}{\sqrt\delta},\sqrt\delta\,x_1\right).
\]
The box ranges imply $c_g\le|h_\delta kx|\le C_g$, also after area normalization. A unimodular lattice containing a primitive vector of length in $[c,C]$ has systole at least $\min(c,C^{-1})$: any independent vector has determinant of magnitude at least one, while a parallel primitive vector agrees up to sign. Mahler compactness therefore gives the fixed compact support.

Choose the bumps by rescaling fixed one-dimensional functions. Their first-box derivatives of order $(j_1,j_2)$ are bounded by
$C_j(\tau\sqrt\delta)^{-j_1}(\sqrt\delta)^{-j_2}$, and derivatives of total order $j$ in the remaining boxes are bounded by $C_j\delta^{-j/2}$. An infinitesimal matrix action differentiates in a direction $Bx$ with coefficients $O_g(\sqrt\delta)$ on these boxes. The one-vector argument applies also to differentiated sums, giving Lie derivatives of $\Phi$ of order at most $j$ bounded by $C_{g,j}\tau^{-j}$. Since $\|\operatorname{Ad}(h_\delta)\|\le C\delta^{-1}$, 
\[
\|\Psi\|_{C^j}\le C_{g,j}\delta^{-j}\tau^{-j}.
\]
On the fixed compact set this implies \eqref{eq:sobolev}. Haar invariance also preserves \eqref{eq:test-mean}. Finally,
\[
kD_Q\ell_{-c}k^{-1}=a_{\frac12\log Q}u_c,
\qquad h_\delta a_{\frac12\log Q}u_c=a_Tu_{c-1/Q},
\]
which proves \eqref{eq:orbit-identity}.
\end{proof}

Put $G=\SL_2(\R)^g$, and let $H<G$ be the diagonal copy of $\SL_2(\R)$. Fix a right-invariant Riemannian metric on $G$, its quotient metric on $X=(\SL_2(\R)/\SL_2(\Z))^g$, and the induced volume normalizations on the orbits. An orbit in $X$ is \emph{periodic} if it is closed and has finite invariant volume.

\begin{theorem}[Effective equidistribution,{\cite[Theorem1.1]{LindenstraussMohammadiYang}}]\label{thm:effective}
There are $A_0,\kappa>0$, a fixed finite-order Sobolev norm $S_0$, and $R_0$, depending on $X$ and the metric, with the following property. Let $z\in X$, let $R_{\log}\ge R_0$ satisfy
\begin{equation}\label{eq:injectivity-condition}
2^{R_{\log}/A_0}>\max\{1,\inj(z),\inj(z)^{-1}\},
\end{equation}
and let $N\in\mathbb N$ with $N\ge A_0R_{\log}$. At least one of the following holds:
\begin{enumerate}
\item For every $f\in C_c^\infty(X)$,
\begin{equation}\label{eq:effective-integer}
\left|\frac12\int_{-1}^1f(a_Nu_bz)\,db-\int_Xf\,dm_X\right|
\le S_0(f)2^{-\kappa R_{\log}}.
\end{equation}
\item There are $x\in X$ and a proper closed subgroup $H\le M<G$ such that $Mx$ is periodic and
\begin{equation*}
\vol(Mx)\le2^{R_{\log}},\qquad
\dist(z,x)\le N^{A_0}2^{-2N+A_0R_{\log}}.
\end{equation*}
\end{enumerate}
\end{theorem}

For real $T$, apply the theorem with $N=\lfloor T\rfloor$ and $f\circ a_{T-N}$. After enlarging $A>A_0$ and the fixed Sobolev norm $S$, the assumptions $T\ge AR_{\log}$ and \eqref{eq:injectivity-condition} with $A$ in place of $A_0$ give either
\begin{equation}\label{eq:effective-real}
\left|\frac12\int_{-1}^1f(a_Tu_bz)\,db-\int_Xf\,dm_X\right|
\le S(f)2^{-\kappa R_{\log}},
\end{equation}
or a periodic orbit of volume at most $2^{R_{\log}}$ with
\begin{equation*}
\dist(z,x)\le T^A2^{-2T+AR_{\log}}.
\end{equation*}
The bounded translations $f\mapsto f\circ a_s$, $0\le s\le1$, are absorbed into $S$. For $\Psi$, their supports remain in a fixed compact set and \eqref{eq:sobolev} holds uniformly. In particular, the exponential factor remains $2^{-2T}$.

\begin{proposition}\label{prop:count}
Suppose \eqref{eq:effective-real} holds for the input \eqref{eq:orbit-input}, and put $R=2^{R_{\log}}$. If
\begin{equation*}
\delta^{-(P+g)}\tau^{-(d_1+1)}R^{-\kappa}=o(1),
\end{equation*}
then the number of children satisfies \eqref{eq:branching}.
\end{proposition}

\begin{proof}
Apply \eqref{eq:effective-real} to $\Psi$ and use \eqref{eq:test-mean}, \eqref{eq:sobolev} and \eqref{eq:orbit-identity}. The error condition gives
\[
\int_I\Phi\bigl(D_Q\ell_{-b}(\Lambda_1,\ldots,\Lambda_g)\bigr)\,db
\ge c_g\delta^g\tau.
\]
By \eqref{eq:Lambda_i} and \eqref{eq:test-functions},
\[
\Phi\bigl(D_Q\ell_{-b}(\Lambda_1,\ldots,\Lambda_g)\bigr)
=\sum_{\mathbf v\in\prod_{i=1}^g L_i^{\mathrm{prim}}}
\prod_{i=1}^g\psi_i(H_{w,b}v_i).
\]
Tonelli's theorem allows termwise integration. For a fixed tuple
with a nonzero contribution, the second target vector has constant
horizontal coordinate at least $\sqrt\delta$ and
$\frac{d}{db}(H_{w,b}v_2)_y=-Q(H_{w,b}v_2)_x$.
Its vertical window has width $2c\sqrt\delta$, so the tuple contributes on a set of parameters of length at most $2c/Q$. Each term is at most one. If $N$ tuples contribute, then $c_g\delta^g\tau\le(2c/Q)N$, and hence
\[
N\ge c_gQ\delta^g\tau=c_g\delta^gW^t.
\]
By the fixed-coordinate window description above and
Lemma~\ref{lem:children}, these tuples give distinct children.
\end{proof}

\section{Normal splittings and periodic obstructions}\label{sec:normal}
Set $\hat{\gamma}=(\gamma-1)/2$ and introduce the earlier expansion $\hat{Q}=W^{\hat{\gamma}}$.
Let
\[
R=W^{\eta \hat{\gamma}},\qquad \sigma=R^{-b_0},\qquad D=R^{A+1},
\]
where $\eta,b_0>0$ will be chosen in \eqref{eq:parameters}.
Increase the external constant $A$ when necessary. Throughout this
section, $z$ is the input in \eqref{eq:orbit-input} and $\Lambda_i$ is in \eqref{eq:Lambda_i}.

\begin{definition}\label{def:normal}
A splitting is \emph{normal} if
$\min_{1\le i\le g}\lambda_1(\Lambda_i)\ge\sigma$ and, for every
pair $i<j$, the comparison
\begin{equation}\label{eq:persistent-comparison}
\lambda_1\bigl((a_{\frac12\log \hat{Q}}u_bz)_j\bigr)
\le C_{\mathrm{per}}D^{5/2}
\lambda_1\bigl((a_{\frac12\log \hat{Q}}u_bz)_i\bigr)
\end{equation}
fails for at least one $b\in[-1,1]$. Here $C_{\mathrm{per}}$ is a
sufficiently large fixed constant.
\end{definition}

We say that a pair satisfies the \emph{persistent comparison} if
\eqref{eq:persistent-comparison} holds for every $b\in[-1,1]$.
The test includes all $\binom g2$ pairs. The next two lemmas explain
why these pair tests detect every proper intermediate orbit.

\begin{lemma}
\label{lem:periodic-pair}
For $G=\SL_2(\R)^g$ and $H<G$ is the diagonal copy of $\SL_2(\R)$,  $M$ is a subgroup that $H\le M<G$ is closed and $Mx$ is periodic. There are indices
$i<j$ such that the diagonal orbit of $(x_i,x_j)$ is periodic and has
volume at most $C_g\vol(Mx)$.
\end{lemma}

\begin{proof}
Write $S=\SL_2(\R)$ and $\Gamma_0=\SL_2(\Z)$.
Every connected intermediate subgroup is a product of diagonal groups
on the blocks of a partition $\mathcal P$ of $\{1,\ldots,g\}$; thus
\[
M^\circ=\prod_{B\in\mathcal P}\Delta_BS,
\qquad
M\subseteq N_G(M^\circ)=M^\circ Z(G).
\]
The center $Z(G)=\{\pm I\}^g$ acts trivially on $X$, so
$Mx=M^\circ x$ with the same induced orbit volume.

Write $x=(g_1\Gamma_0,\ldots,g_g\Gamma_0)$. The stabilizer in the
block $B$ is
$\Gamma_B=\bigcap_{i\in B}g_i\Gamma_0g_i^{-1}$.
Periodicity makes every $\Gamma_B$ a lattice, and the fixed volume
normalizations give
\[
\vol(Mx)\asymp_g\prod_{B\in\mathcal P}\covol_S(\Gamma_B),
\qquad
\covol_S(\Gamma_B)\ge\covol_S(\Gamma_0)>0.
\]
Since $M$ is proper, some block contains distinct indices $i,j$. For
this block,
\[
\Gamma_B\subseteq
\Gamma_{ij}:=g_i\Gamma_0g_i^{-1}\cap g_j\Gamma_0g_j^{-1},
\qquad
\covol_S(\Gamma_{ij})\le\covol_S(\Gamma_B).
\]
The lower bounds for the other block covolumes therefore bound the
pair-orbit volume by $C_g\vol(Mx)$. Moreover, $\Gamma_{ij}$ has
finite index in $g_i\Gamma_0g_i^{-1}$. The first-coordinate projection
of $S/\Gamma_{ij}$ is a finite covering of $X_1$, so the pair-orbit
map is proper. Its image is consequently closed and has the finite
volume just bounded.
\end{proof}

\begin{lemma}\label{lem:commensurable}
If a diagonal orbit $\mathcal O\subset X_1^2$ is periodic and
$\vol(\mathcal O)\le L$, then every $(\Lambda_i,\Lambda_j)$ on
that orbit satisfies
\begin{equation}\label{eq:periodic-systoles}
(C\sqrt L)^{-1}\lambda_1(\Lambda_i)
\le\lambda_1(\Lambda_j)
\le C\sqrt L\,\lambda_1(\Lambda_i).
\end{equation}
The constant depends only on the fixed volume normalization.
\end{lemma}

\begin{proof}
Write the pair as $(h\Z^2,hz'\Z^2)$ and put $\Gamma_0=\SL_2(\Z)$.
Periodicity implies that $z'$ commensurates $\Gamma_0$.
The standard commensurator description and Hecke index bound give
$z'=d^{-1/2}B$ where
$B\in \operatorname{Mat}_{2\times2}(\mathbb Z)$ is primitive and 
$\det B=d>0$. $d$ also satisfies
\[
d\le[\Gamma_0:\Gamma_0\cap z'\Gamma_0(z')^{-1}]
\ll\vol(\mathcal O)\le L.
\]
The inclusions $d\Z^2\subseteq B\Z^2\subseteq\Z^2$ yield
\[
\lambda_1(h\Z^2)\le\lambda_1(hB\Z^2)
\le d\,\lambda_1(h\Z^2).
\]
Multiplication by $d^{-1/2}$ proves \eqref{eq:periodic-systoles},
uniformly in $h$.
\end{proof}

\begin{lemma}\label{lem:obstruction}
With the parameters chosen as in \eqref{eq:parameters}, every
sufficiently tall normal input satisfies \eqref{eq:effective-real}.
\end{lemma}

\begin{proof}
Standard reduction theory for unimodular lattices, together with the
bounded shear and rotation defining $z$, gives
\begin{equation*}
c_g\min\left\{1,\min_i\lambda_1(\Lambda_i)^2\right\}
\le\inj(z)\le C_g.
\end{equation*}
Normality and $2b_0<1/A$ imply
$\inj(z)^{-1}\le C R^{2b_0}=o(R^{1/A})$.
The parameter choices also give $T\ge A R_{\log}$ for large $W$.
Thus the size hypotheses of Theorem~\ref{thm:effective} hold.

Suppose its periodic alternative occurs, and put
$\rho=T^A2^{-2T}R^A$. We have $\rho\to0$, since
$T=\hat{\gamma}\log W+O(\log\log W)$ and $A\eta<\log 2/4$.
The quotient metric gives representatives with $z=hx$ and
\[
\|h-I\|+\|h^{-1}-I\|\le C_g\rho.
\]
For $|b|\le1$, conjugation by $a_{\frac12\log \hat{Q}}$ amplifies
a perturbation by at most $\hat{Q}$. Hence
\begin{equation}\label{eq:perturbation-bound}
\left\|
\left(a_{\frac12\log \hat{Q}}u_bh u_{-b}
a_{-\frac12\log \hat{Q}}\right)^{\!\pm1}
\right\|
\le C_g(1+\hat{Q}\rho)\le D.
\end{equation}
For the last inequality, use $\hat{Q}\le e^T$ and $2\log2>1$ to get
$\hat{Q}\rho\le T^AR^A$, while $T^A=o(R)$.

The unperturbed arc remains in $Mx$, since its translating matrices
belong to $H\le M$. Lemma~\ref{lem:periodic-pair} supplies one fixed
pair $i<j$ on a periodic pair orbit of volume at most $C_gR$.
Lemma~\ref{lem:commensurable} compares its systoles by $C_g\sqrt R$
at every time. If $h_{b,i},h_{b,j}$ are the corresponding displacement
factors from \eqref{eq:perturbation-bound}, then
\begin{equation}\label{eq:perturbed-systoles}
\lambda_1(h_{b,j}\Lambda_j')
\le D\lambda_1(\Lambda_j')
\le C_gD\sqrt R\,\lambda_1(\Lambda_i')
\le C_gD^2\sqrt R\,\lambda_1(h_{b,i}\Lambda_i').
\end{equation}
Since $R\le D$, the comparison factor is at most
$C_{\mathrm{per}}D^{5/2}$. Thus this same pair satisfies
\eqref{eq:persistent-comparison} for every $b\in[-1,1]$, contrary
to normality. This excludes the periodic alternative.
\end{proof}

\section{Arrival measures and period fibers}\label{sec:fibers}
Throughout this section, $g\ge3$ is fixed and $r\varepsilon<t$. We use the
parameter choices in \eqref{eq:parameters} and take $W_0$ beyond the threshold
verified in Section~\ref{sec:survival}. At a splitting with slit $w$ and
actual height $W=|w|$, recall the normalization \eqref{eq:Lambda_i}:
\begin{equation*}
 G_w=\diag\left(
 \sqrt{\frac{d}{W^{1+\gamma}\log W}},
 \sqrt{\frac{W^{1+\gamma}\log W}{d}}\right)R_{-\arg w},
 \qquad \Lambda_i=A_i^{-1/2}G_wL_i.
\end{equation*}
The children in $\mathcal C(S)$ of a normal splitting satisfies
Lemma~\ref{lem:children} and Proposition~\ref{prop:count}; in particular,
\begin{equation*}
 b(S):=\#\mathcal C(S)\ge c\delta(W)^gW^t\ge1.
\end{equation*}
A path chooses each candidate with probability $1/b(S)$ and stops at its first
nonnormal input. An integer history $h$ records the primitive loop labels at
each step in the initial marked absolute homology. The update
\eqref{eq:move} determines all later homology lattices and relative slit labels.
Let $D_h$ be the root domain on which these choices satisfy the geometric
windows and all strict ancestors are normal. The endpoint is untested, and its
incoming weight is
\begin{equation}\label{eq:history-weights}
 p_h(Z)=\mathbf 1_{D_h}(Z)\prod_{j<n}b(S_j^h(Z))^{-1}
 \le C^n\prod_{j<n}\delta_j^{-g}W_j^{-t},\qquad |h|=n.
\end{equation}

Choose a marked root region $\Omega$ with normalized lattices in fixed compact
charts, component areas $A_i\ge2a$, and $W_0\le|w|<2W_0$, where
$0<a<1/(2g)$. Let $\mu_\Omega$ be its normalized area-one period measure.
Write $\mathcal S$ for the space of marked cyclic splittings,
$\mathcal N$ for the normal inputs, and $\mathcal B=\mathcal S\setminus\mathcal N$.
Define the arrival measures by
\begin{equation*}
 \nu_{n,Z}=\sum_{\substack{|h|=n\\Z\in D_h}}p_h(Z)\delta_{S_h(Z)},
 \qquad
 \nu_n(E)=\int_\Omega\nu_{n,Z}(E)\,d\mu_\Omega(Z),
 \qquad \nu_{0,Z}=\delta_Z.
\end{equation*}
Distinct candidates give distinct child splittings, and hence
\begin{equation*}
 \nu_{n+1}(E)=\int_{\mathcal N}
 \frac{\#(\mathcal C(S)\cap E)}{b(S)}\,d\nu_n(S).
\end{equation*}
The same recursion holds for $\nu_{n,Z}$. These measures are Borel: integer
loop labels are countable, membership in the open orbit windows can be tested
at rational times, and continuity of the systole reduces the persistent
comparison to rational times. The finite candidate counts, history domains,
and weights are therefore Borel.

Put $\ell_n(Z)=\nu_{n,Z}(\mathcal B)$ and $L_n=\nu_n(\mathcal B)$.
The mass surviving through depth $N$ is
$s_N(Z)=\nu_{N,Z}(\mathcal N)=1-\sum_{n=0}^N\ell_n(Z)$.
Thus $s_N\downarrow s$ for a Borel function $s\in[0,1]$, and Tonelli gives
\begin{equation}\label{eq:conservation}
 \int_\Omega s(Z)\,d\mu_\Omega(Z)
 =1-\sum_{n\ge0}\nu_n(\mathcal B).
\end{equation}
All weights remain the original incoming weights; discarded mass is not
redistributed among the surviving children.

For the period-coordinate estimates, set
\begin{equation*}
 K_0=W_0^{1+\gamma}\log W_0,\qquad
 G_*=\diag(K_0^{-1/2},K_0^{1/2})R_*,\qquad
 Z_*=(P_*,z_*)=(G_*P,G_*w),
\end{equation*}
where $P\in\R^{2\times2g}$ consists of the marked absolute periods. Divide the
root directions into at most $CK_0$ sectors of aperture $c/K_0$ and let $R_*$
send the central direction to the horizontal. On each sector,
$\|P_*\|\le C$ and $|z_*|\le CW_0/\sqrt{K_0}$.
We refer to the coordinates obtained by applying $G_*$
to all period vectors as the \emph{frozen frame}.
The unnormalized measure of
$\Omega$ is comparable to $W_0^2$: the area and normalized lattice parameters
range over fixed compact sets, while the slit ranges over an annulus of that
area. At a fixed slit the common map $G_w$ has determinant one on each period
vector, and its dependence on the slit gives a block triangular change of
variables with determinant one. Thus normalization contributes a factor
comparable to $W_0^{-2}$.

Partition histories into compatible dyadic height words
$\mathbf W=(W_0,\ldots,W_n)$, where $W_j\le|w_j|<2W_j$.
There are at most $C^n$ such words since
$|\log W_{j+1}-\gamma\log W_j|\le C$.
Let $\nu_{n,\mathbf W}$ be the corresponding restricted arrival measure.
Until stated otherwise, fix one word and one root sector, and restrict all
history domains and weight suprema accordingly. The height recurrence gives
\[
 \log|w_j|=\gamma^j\log|w_0|
 +O\!\left(\frac{\gamma^j-1}{\gamma-1}\right),
 \qquad
 \sum_j\delta_j\le\frac{C}{(\gamma-1)\log W_0}.
\]
For sufficiently large $W_0$, the total area exchange is less than $a$, so
every component area stays in $[a,1-(g-1)a]$. The total turn is $O(K_0^{-1})$;
all subsequent slits and chosen loops remain in a fixed narrow root cone, on
which $|G_*q|\asymp|q|/\sqrt{K_0}$.

Fix the integral symplectic basis of absolute homology
determined by the initial marked splitting, and identify the resulting coefficient lattice with $\Z^{2g}$.
For each primitive loop with holonomy $v_i$ selected on the last edge of a history, let $h_i\in\Z^{2g}$ denote the coordinate vector of its homology class in this fixed basis.
This loop is retained in the $i$th child torus. 
Because $P$ is the matrix of marked absolute periods, then
\begin{equation*}
    v_i=Ph_i,\qquad q_i=P_*h_i=G_*v_i.
\end{equation*}
For a fixed integer history, the labels $h_i$ remain fixed as the period coordinates vary.
For the
$i$-th child torus with lattice $L_i$,
complete $h_i$ to a positively oriented integral basis
$(h_i,k_i)$ of $L_i$, and define
$\Delta_i:=h_i\wedge k_i\in\bigwedge^2\Z^{2g}$.
This bivector is independent of the chosen oriented basis,
and $|\Delta_i|$ is the Euclidean covolume of $L_i$
in its real span.
\begin{lemma}\label{lem:history}
Fix a height word of depth $n\ge1$. On its last edge write
\[
 W_-=W_{n-1},\qquad V=W_-^r,\qquad W=W_n\asymp W_-^\gamma,
 \qquad M=V/\sqrt{K_0},\qquad N=W/\sqrt{K_0}.
\]
Then 
\begin{equation}\label{eq:label-sizes}
 \begin{aligned}
 |h_1|&\asymp M,& |h_i|&\asymp N\quad(2\le i\le g),&
 |h_1\wedge\cdots\wedge h_g|&\asymp MN^{g-1},\\
 |\Delta_1|&\asymp MN,& |\Delta_i|&\asymp N^2\quad(2\le i\le g).
 \end{aligned}
\end{equation}
The child homology lattices have bases of label length at most $CN$.
Writing the slit label as $(m,1)$ gives
$m=m_-+\sum_i h_i$, with $|m|\le CN$ and
$|m_-|\le CW_-/\sqrt{K_0}=o(M)$. Moreover, for
$E_n=\exp(Cn^2+Cn\log\log W_0+Cn)$,
\begin{equation}\label{eq:weighted-histories}
 \Xi_n:=\sum_{|h|=n}\sup_{D_h}p_h
 \le CM^2N^{2g-2}\frac{W_0}{V}E_n.
\end{equation}
Empty history domains contribute zero.
\end{lemma}

\begin{proof}
Fix the height word and root sector, and write $K=K_0$.
Put $V_j=W_j^r$, $U_j=W_j^\gamma$, and
$\delta_j=d/\log W_j$.
Let $v_{i,j}$ be the loop selected at depth $j$, with
homology label $h_{i,j}$, and let
$\mathcal L_i^{(j)}\subset\Z^{2g}$ be the $i$th homology
lattice at that depth.
Thus $h_i=h_{i,n-1}$, while the bivector $\Delta_i$ in
the statement belongs to $\mathcal L_i^{(n)}$.
All estimates below are uniform over the corresponding
history domains.

Let $P:\R^{2g}\to\R^2$ be the absolute period map.
Its restriction to the real span of $\mathcal L_i^{(j)}$
is an isomorphism.
Denote its inverse by
$F_{i,j}:\R^2\to\R^{2g}$, so that
$F_{i,j}v_{i,j}=h_{i,j}$.
If $B_i$ is the initial $i$th period matrix in the frozen
frame and $\iota_i$ is the inclusion into the $i$th
coordinate block, then
$F_{i,0}=\iota_iB_i^{-1}G_*$.

\eqref{eq:move} gives the exact identity
\[
F_{i,j+1}x-F_{i,j}x
=
\frac{v_{i,j}\times x}{A_i^{(j+1)}}
\sum_{k\ne i}(F_{i,j}-F_{k,j})v_{k,j},
\qquad x\in\R^2.
\]
Indeed, both sides vanish on $v_{i,j}$ and agree on the
new complementary period
$u_{i,j}-\sum_{k\ne i}v_{k,j}$.

Let $I_j\subset S^1$ be the angular interval centered on
$w_j/|w_j|$, of radius
$r_j=C_{\mathrm{ang}}\delta_j/(W_jU_j)$.
Choose $C_{\mathrm{ang}}$ sufficiently large.
The angle estimates in Lemma~\ref{lem:children} place every
$v_{i,j}/|v_{i,j}|$ in $I_j$.
Moreover, the slit-turn estimate and
$r_{j+1}/r_j\to0$ imply
$I_{j+1}\subset I_j$ for sufficiently large $W_0$.
For $e\in I_{j+1}$, the child determinant estimate gives
\[
|v_{i,j}\times e|
\le
\frac{|v_{i,j}\times w_{j+1}|}{|w_{j+1}|}
+
|v_{i,j}|r_{j+1}
\le C\frac{\delta_j}{U_j}.
\]

Set
$C_j=\sqrt K\max_i\sup_{e\in I_j}|F_{i,j}e|$.
The compact initial charts and the frozen-cone estimate
give $C_0\le C$.
Using $A_i^{(j+1)}\ge a$ and
$\sum_k|v_{k,j}|\le CU_j$ in the update identity yields
$C_{j+1}\le(1+C\delta_j)C_j$.
Since $\sum_j\delta_j$ is uniformly small, we have
$C_j\le C$.
The nesting of the intervals then allows us to telescope
at a fixed direction:
\[
|(F_{i,j}-F_{i,0})e|
\le \frac C{\sqrt K}\sum_{k<j}\delta_k,
\qquad e\in I_j.
\]
Apply this to the direction of the last-edge loop $v_i$.
Since $|G_*v_i|\asymp |v_i|/\sqrt K$, we obtain
\begin{equation}\label{eq:block-approximation}
\frac{h_i}{|G_*v_i|}
=
\iota_i(B_i^{-1}e_i)
+
O\!\left(\sum_{j<n}\delta_j\right),
\qquad
e_i=\frac{G_*v_i}{|G_*v_i|}.
\end{equation}
The same estimate holds at every depth.
The leading vectors occupy mutually orthogonal coordinate
blocks, and their lengths are uniformly bounded above and
below.
Taking $W_0$ sufficiently large therefore gives
$|h_1|\asymp M$, $|h_i|\asymp N$ for $i\ge2$, and
$|h_1\wedge\cdots\wedge h_g|\asymp MN^{g-1}$.

Write $\Delta_i^{(j)}$ for the oriented bivector of
$\mathcal L_i^{(j)}$.
The basis update gives
\[
\Delta_i^{(j+1)}
=
\Delta_i^{(j)}-E_{i,j},
\qquad
E_{i,j}:=
h_{i,j}\wedge\sum_{k\ne i}h_{k,j}.
\]
The block approximation implies uniform independence of
the normalized labels. Consequently,
\[
|E_{i,j}|
\asymp
|h_{i,j}|
\left(\sum_{k\ne i}|h_{k,j}|^2\right)^{1/2},
\]
and hence
\[
|E_{1,j}|\asymp\frac{V_jU_j}{K},
\qquad
|E_{i,j}|\asymp\frac{U_j^2}{K}
\quad(i\ge2).
\]
Here $g\ge3$ ensures that, for each $i\ge2$, the sum
contains another long label.

The ratios of consecutive increment scales tend uniformly
to infinity as $W_0\to\infty$.
Their initial values are
$W_0^t/\log W_0$ and
$W_0^{\gamma-1}/\log W_0$, respectively, and also tend to
infinity.
Thus a geometric-series estimate gives, uniformly in $n$,
\[
|\Delta_i^{(0)}|
+
\sum_{j<n-1}|E_{i,j}|
=
o\bigl(|E_{i,n-1}|\bigr).
\]
The last increment therefore determines the size of the
child bivector.
Using $W_n\asymp U_{n-1}$, we obtain
$|\Delta_1|\asymp MN$ and
$|\Delta_i|\asymp N^2$ for $i\ge2$.

The retained label $h_i$ is primitive in the child lattice.
Complete it to a basis $(h_i,k_i)$ and reduce $k_i$ modulo
$\Z h_i$ so that its component parallel to $h_i$ has
length at most $|h_i|/2$.
Its perpendicular component has length
$|\Delta_i|/|h_i|$, and therefore
\[
|k_i|
\le \frac{|h_i|}{2}+\frac{|\Delta_i|}{|h_i|}
\le CN.
\]
This proves the basis bound.

Finally, $m_0=0$ and
$m_{j+1}=m_j+\sum_i h_{i,j}$.
Since $W_{j+1}\asymp U_j$ and the heights grow at least
geometrically,
\[
|m_n|
\le
\frac C{\sqrt K}\sum_{j<n}U_j
\le \frac{CW_n}{\sqrt K}=CN.
\]
Similarly,
$|m_{n-1}|\le CW_{n-1}/\sqrt K=o(M)$, proving the
relative-label assertions.

A rank-two Euclidean lattice of covolume $D$ and with
a basis bounded by $B$ has $O(R^2/D)$ points in a disk
of radius $R\ge2B$, as follows by translating a
fundamental parallelogram.

Set
$\mathcal F_j:=V_j^2W_{j+1}^{2g-2}/K^g$.
At the root, the number of possible loop tuples is at most
$C\mathcal F_0$.
For $j\ge1$, the current homology covolumes are comparable
to $V_{j-1}W_j/K$ and $W_j^2/K$.
The preceding lattice-point estimate therefore bounds
the number of extensions of any fixed prefix by
\[
C\frac{V_j^2U_j^{2g-2}}
        {V_{j-1}W_j^{2g-1}}
\le
C\frac{\mathcal F_j}{\mathcal F_{j-1}}
 \frac{V_{j-1}}{W_j}.
\]
The counting radii dominate the old basis bounds because
$V_j/W_j=W_j^t\to\infty$.
At the root this follows from
$V_0^2/K=W_0^{t-r\varepsilon}/\log W_0\to\infty$.

For $0\le j\le n$, let $\mathscr H_j$ be the set of
admissible integer histories of length $j$, restricted
to the fixed root sector and the height word
$(W_0,\ldots,W_j)$. Define
$\Xi_j=\sum_{h\in\mathscr H_j}\sup_{Z\in D_h}p_h(Z)$,
where $p_h$ is the incoming weight of the prefix $h$.
Each distinct prefix is counted once, and $\Xi_0=1$.
A child domain is contained in its prefix domain, and
the branching bound multiplies its weight supremum by
at most $C\delta_j^{-g}W_j^{-t}$.
Since $r=1+t$, we obtain
\[
\frac{V_0\Xi_1}{\mathcal F_0}
\le CW_0\delta_0^{-g},
\qquad
\frac{V_j\Xi_{j+1}}{\mathcal F_j}
\le
C\delta_j^{-g}
\frac{V_{j-1}\Xi_j}{\mathcal F_{j-1}}
\quad(j\ge1).
\]
Iterating and using
$\mathcal F_{n-1}=M^2N^{2g-2}$ gives
\[
\Xi_n
\le
C^nM^2N^{2g-2}\frac{W_0}{V}
\prod_{j<n}\delta_j^{-g}.
\]
Finally,
$\log W_j\asymp\gamma^j\log W_0$ implies
$\sum_{j<n}\log\delta_j^{-g}
\le Cn^2+Cn\log\log W_0+Cn$.
Absorbing the remaining factors $C^n$ into $E_n$
proves \eqref{eq:weighted-histories}.
\end{proof}

\begin{lemma}\label{lem:fiber}
For each history in a fixed height word, let $F_h$ be an endpoint event in its
original domain. Fix the frozen outputs $y=Z_*(m,1)$ and $q_i=P_*h_i$ for $1\le i\le g$.
Choose fixed symplectic complements $k_i$ in the child tori and write
\begin{equation}\label{eq:literal-projections}
 q_{k_i}=\xi_iq_i+A_i\frac{Jq_i}{|q_i|^2}.
\end{equation}
The $\xi_i$ are literal projections, without reduction modulo integers.
Suppose that, on a fixed bounded radial thickening of the root domain, the
event has a containing extension whose allowed area vectors
$(A_1,\ldots,A_g)$ have $g$-dimensional volume at most $\kappa$ for every fixed
output and every fixed $(\xi_1,\ldots,\xi_g)$. Then
\begin{equation}\label{eq:fiber-bound}
 \int_\Omega\sum_{|h|=n}p_h\mathbf1_{F_h}\,d\mu_\Omega
 \le C\delta_{n-1}^g\kappa E_n(1+W_0^{r\varepsilon}).
\end{equation}
\end{lemma}

\begin{proof}
The retained loops are disjoint, so their independent homology classes span
a Lagrangian plane $E\subset\R^{2g}$. Choose complements with
$\langle h_i,k_j\rangle=\delta_{ij}$ and put $e_i=\Pi_{E^\perp}k_i$.
The standard symplectic matrix identifies this duality with Euclidean
duality between $E$ and $E^\perp$. The Gram bound in
\eqref{eq:block-approximation} therefore gives
\begin{equation*}
 |e_1|\le C/M,\qquad |e_i|\le C/N\quad(i\ge2).
\end{equation*}

We first control the relative label's distance from $E$.
Write a next loop as $h_i'=\alpha_i h_i+\beta_i k_i$. Physical determinants
give $\beta_i=(v_i\times v_i')/A_i$. Relative to the current slit of height
$W$, the retained loop has transverse component $O(\delta/W)$, with
$\delta=\delta_{n-1}$; the next loop has an at most comparable transverse
component and greater length. Hence $|\beta_i|\le C\delta|v_i'|/W$.
Its perpendicular label component is $\beta_i e_i$, so after division by
$|h_i'|\asymp|v_i'|/\sqrt{K_0}$ the bounds are $C\delta/(MN)$ for $i=1$
and $C\delta/N^2$ for $i\ge2$. The uniform Gram bound for the next labels
then implies
$\|\Pi_{E'}-\Pi_E\|\le C\delta/(MN)$.

Index the retained plane and relative label after $j$ moves by $E_j,m_j$,
and put $N_j=W_j/\sqrt{K_0}$ and $M_j=V_{j-1}/\sqrt{K_0}$.
Since $m_{j+1}-m_j\in E_{j+1}$ and $|m_j|\le CN_j$,
\begin{equation*}
 \dist(m_{j+1},E_{j+1})
 \le\dist(m_j,E_j)+C\frac{\delta_{j-1}}{M_j}.
\end{equation*}
Starting from $m_1=\sum_i h_i^1\in E_1$ and using the geometric growth of
$V_j$, we obtain
\begin{equation*}
\dist(m,E)\le\rho_0:=C\sqrt{K_0}\sum_{j\ge0}\frac{\delta_j}{V_j}
 \le C\delta_0\frac{\sqrt{K_0}}{V_0}=o(1).
\end{equation*}
The last limit uses $r\varepsilon<t$. Write
$m=\sum_i a_i h_i+r_\perp$, where $r_\perp\perp E$.
Since $m=m_-+\sum_i h_i$ and $|m_-|\le CW_-/\sqrt{K_0}$, the same dual
bounds imply $|a_i-1|\le CW_-/V=o(1)$.

The pairs $(h_i,k_i)$ form an integral symplectic basis, so the change to
their periods has determinant of absolute value one. The coefficient of the
initial relative period in $y$ is one, and, for fixed $q_i$, the Jacobian
of $(\xi_i,A_i)\mapsto q_{k_i}$ in \eqref{eq:literal-projections} is one.
Thus full-period Lebesgue measure is exactly
\begin{equation}\label{eq:fiber-jacobian}
 dZ_*=dy\prod_{i=1}^g dq_i\,d\xi_i\,dA_i,
\end{equation}
where $dy$ and $dq_i$ are planar Lebesgue measures.

Put $q_{k_i}:=P_*k_i$ and  
write $k_i=\sum_{j=1}^g c_{ij}h_j+e_i$,
where the coefficients $c_{ij}$ depend only on the
fixed homology labels and the chosen complements.
Applying $P_*$ gives
$q_{k_i}=\sum_{j=1}^g c_{ij}q_j+P_*e_i$.
Since $\|P_*\|\le C$,
\begin{equation}\label{eq:projection-intervals}
 \left|\xi_i-
 \frac{q_i\cdot\sum_jc_{ij}q_j}{|q_i|^2}\right|
 \le C\frac{|e_i|}{|q_i|}
 \le
 \begin{cases}
 CM^{-2},&i=1,\\
 CN^{-2},&2\le i\le g.
 \end{cases}
\end{equation}
These interval centers depend only on the outputs and labels, not on the
areas. Enlarge the domain to the product of these intervals and a fixed
compact box of positive areas. The assumed area estimate then bounds its
$(\xi,A)$-volume by $C\kappa M^{-2}N^{-2g+2}$.

Set $e_0=C(W_0/\sqrt{K_0}+\rho_0)$. The relative decomposition gives
$|y-\sum_i a_iq_i|=|z_*+P_*r_\perp|\le e_0$.
Successively fixing the outputs, their containing areas are
\begin{equation}\label{eq:output-areas}
 \begin{aligned}
 y&:\ CN^2,& q_1&:\ C\delta M/N,\\
 q_i&:\ C\delta\quad(2\le i\le g-1),&
 q_g&:\ Ce_0\delta/N.
 \end{aligned}
\end{equation}
Indeed, $|y\times q_i|\le C\delta$ gives transverse width $C\delta/N$,
while the longitudinal ranges are $O(M)$ and $O(N)$. Once
$y,q_1,\ldots,q_{g-1}$ are fixed, $|a_g|\ge c$ confines $q_g$ to a disk
of radius $Ce_0$, giving the final saving. The product of the output areas
is $C\delta^g e_0M$. Together with \eqref{eq:fiber-jacobian}, this bounds
one history's full-period volume by
$C\kappa\delta^g e_0/(MN^{2g-2})$.

For an area-one set $F$, let
$\widehat F=\{\rho Z:Z\in F,\ 1\le\rho\le1+\kappa_0\}$, where
$\kappa_0>0$ is fixed. Radial disintegration gives
$\mu_\Omega(F)\le CW_0^{-2}\vol_{4g+2}(\widehat F)$.
Dilation preserves the labels and changes the containing geometric bounds
only by fixed factors. Keep each original weight supremum fixed on this cone.
Multiplying the preceding volume by \eqref{eq:weighted-histories}, by the
normalizing density $C/W_0^2$, and by the number $CK_0$ of sectors leaves
the scale factor
\begin{equation*}
 \frac{e_0}{MN^{2g-2}}\,
 M^2N^{2g-2}\frac{W_0}{V}\frac{K_0}{W_0^2}
 =\frac{e_0\sqrt{K_0}}{W_0}
 \le C(1+W_0^{\gamma-r})=C(1+W_0^{r\varepsilon}).
\end{equation*}
This proves \eqref{eq:fiber-bound}.
\end{proof}

\begin{lemma}\label{lem:area}
Fix a compact interval $I\subset(0,\infty)$ and $c_*>0$. There are constants
$H_0,C>0$, depending only on $I,c_*$, with the following property.
Suppose $a,b\in\R^2$, $a\times b=1$, $H=|a|\ge H_0$, and
$|c|\ge c_*$, where $c=(a\cdot b)/H^2$. Then, for every $\xi\in\R$ and
every $\sigma>0$,
\begin{equation}\label{eq:area-short-vector}
 \operatorname{Leb}\left\{A\in I:
 \lambda_1\bigl(A^{-1/2}[\Z a+\Z(\xi a+Ab)]\bigr)<\sigma\right\}
 \le C\sigma^2.
\end{equation}
The same bound holds with a non-strict inequality.
\end{lemma}

\begin{proof}
It suffices to consider $0<\sigma\le1$. Rotate so that
$a=(H,0)$ and $b=(Hc,H^{-1})$. A lattice vector with integer coefficients
$(\alpha,\beta)\ne(0,0)$ is
\[
 \left(A^{-1/2}H[\alpha+\beta(\xi+Ac)],\ A^{1/2}\beta/H\right).
\]
Choose $H_0/\sqrt{\max I}>1$. A vector shorter than $\sigma$ then satisfies
\begin{equation*}
 0<|\beta|\le C_IH\sigma,\qquad
 |\alpha+\beta(\xi+Ac)|\le C_I\sigma/H.
\end{equation*}
For each nonzero $\beta$, at most
$C_I(1+|\beta c|)\le C_{I,c_*}|\beta c|$ integers $\alpha$ occur as $A$
ranges over $I$. Each restricts $A$ to an interval of length at most
$C_I\sigma/(H|\beta c|)$. The allowed length for one $\beta$ is therefore
at most $C\sigma/H$. There are at most $CH\sigma$ possible nonzero
$\beta$, proving the bound uniformly in $\xi$. For the non-strict version,
apply the bound at thresholds $\sigma+\epsilon$ and let $\epsilon$ decrease
to zero.
\end{proof}

\begin{proposition}\label{prop:thin}
For a fixed height word and $0<\sigma_n\le1$, the original incoming mass
$L_n^{\mathrm{thin}}(\mathbf W;\sigma_n)$ of
$\min_i\lambda_1(\Lambda_i)<\sigma_n$ satisfies
\begin{equation}\label{eq:thin-loss}
 L_n^{\mathrm{thin}}(\mathbf W;\sigma_n)
 \le C\delta_{n-1}^g E_n(1+W_0^{r\varepsilon})\sigma_n^2.
\end{equation}
\end{proposition}

\begin{proof}
At fixed outputs put $T=G_wG_*^{-1}$ and, for the $i$th loop, set
$a=Tq_i$, $b=TJq_i/|q_i|^2$, $H=|a|$, and $c=a\cdot b/H^2$.
These vectors are fixed while the areas vary, and $a\times b=1$.
Use $L=V$ for $i=1$ and $L=W$ for $i\ge2$. All estimates use the actual
current height $W$, comparable to its word representative, and put
$\delta=\delta_{n-1}\asymp1/\log W$. Write \(R_{-\arg w}v_i=(v_{i,\parallel},v_{i,\perp})\), so $
v_i=v_{i,\parallel}\frac{w}{W}+v_{i,\perp}\frac{Jw}{W}.
$
The incoming signed windows give $|v_{i,\parallel}|\asymp L$ and
$|v_{i,\perp}|\asymp\delta/W$. Thus
\[
 a_\parallel=v_{i,\parallel}
 \sqrt{\frac d{W^{1+\gamma}\log W}},\qquad
 a_\perp=v_{i,\perp}
 \sqrt{\frac{W^{1+\gamma}\log W}{d}},
 \qquad
 \frac{|a_\parallel|}{|a_\perp|}
 \le CW^{1-\gamma}\longrightarrow0.
\]
We claim that
\begin{equation*}
 H^2\asymp\delta^2W^{\gamma-1}\log W
 \asymp\frac{W^{\gamma-1}}{\log W}\longrightarrow\infty,
 \qquad |c|\asymp\frac W{L\delta}\longrightarrow\infty.
\end{equation*}
The first comparison follows from the preceding
coordinate estimates. To prove the second, recall that $w$ is the current slit, 
and $y=G_*w$, while $q_i=G_*v_i$. Then we have \[
q_i=\frac{v_{i,\parallel}}{W}y+v_{i,\perp}G_*\!\left(\frac{Jw}{W}\right).
\]
For the second comparison, the frozen relation
$q_i=(v_{i,\parallel}/W)y+O(\delta\sqrt{K_0}/W)$
has relative error $O(\delta K_0/(WL))=o(1)$, since
$W\asymp W_-^\gamma$, $L\ge cW_-^r$, and $W_-\ge W_0$.
Hence $y\cdot q_i\asymp WL/K_0$ and $|q_i|\asymp L/\sqrt{K_0}$.
Using $\det(y,Jq_i)=y\cdot q_i$ gives
$|b_\perp|\asymp\sqrt{W^{1+\gamma}\log W}/L$.
It follows that
\[
 |a_\perp b_\perp|\asymp\delta W^\gamma\log W/L\longrightarrow\infty,
 \qquad
 |a_\parallel b_\parallel|
 \le CK_0/(W^{1+\gamma}\log W)=O(1).
\]
Dividing the dominant transverse product by $H^2$ proves the claim for $c$.

At fixed literal projections, the $i$th thin condition restricts only $A_i$
to a set of length $C\sigma_n^2$ by Lemma~\ref{lem:area}.
A union bound over the coordinates in a fixed compact area box therefore
gives allowed area volume at most $C\sigma_n^2$. Under a bounded dilation
$\rho$, each normalized lattice changes by the bounded common matrix
$\diag(\kappa_\rho^{-1},\kappa_\rho)$, where
\begin{equation*}
 \kappa_\rho=\rho^{(1+\gamma)/2}
 \sqrt{\frac{\log(\rho W)}{\log W}}\asymp1.
\end{equation*}
Thus the radial cone is contained in the event with threshold $C\sigma_n$.
Lemma~\ref{lem:area} allows every positive threshold, so
Lemma~\ref{lem:fiber} applies with $\kappa=C\sigma_n^2$.
\end{proof}

\begin{lemma}\label{lem:two-area}
Let $\hat{Q}(W)=W^{\hat{\gamma}}$ with fixed $\hat \gamma>0$, and put
$F_b=a_{\frac12\log \hat{Q}}u_bk$, where $b$ ranges over a fixed
bounded interval. For distinct $i,j$, at fixed outputs, literal projections,
and time, the allowed area vectors have volume at most
$C\epsilon_i^2\epsilon_j^2$ on the event
\begin{equation*}
 \lambda_1((F_b\Lambda)_i)\le\epsilon_i,\qquad
 \lambda_1((F_b\Lambda)_j)\le\epsilon_j.
\end{equation*}
The constant is uniform in the outputs, time, history, and positive
thresholds, with the parameters and genus fixed.
\end{lemma}

\begin{proof}
Use the fixed-output basis from the preceding proof, with first vector
$a=(a_\parallel,a_\perp)$ and complementary vector
$d=ca+Ja/H^2$. Since $|a_\parallel|/|a_\perp|\to0$ and $k=J$,
\begin{equation*}
 F_ba=\left(\sqrt{\hat{Q}}(-a_\perp+ba_\parallel),
 \frac{a_\parallel}{\sqrt{\hat{Q}}}\right),
 \qquad |F_ba|\asymp\sqrt{\hat{Q}}\,H.
\end{equation*}
Also $|F_bJa|\le C\sqrt{\hat{Q}}\,H$, so the new shear coefficient is
\begin{equation*}
 c_b=\frac{F_ba\cdot F_bd}{|F_ba|^2}
 =c+\frac{F_ba\cdot F_bJa}{H^2|F_ba|^2}
 =c+O(H^{-2}).
\end{equation*}
Thus Lemma~\ref{lem:area} applies uniformly throughout the bounded time
interval. At fixed outputs, projections, and time, the lattices in factors
$i$ and $j$ depend only on $A_i$ and $A_j$, respectively. Their allowed
lengths are bounded by $C\epsilon_i^2$ and $C\epsilon_j^2$; the remaining
areas range over fixed compact intervals. Fubini gives the product bound in
full period space, as required by Lemma~\ref{lem:fiber}.
\end{proof}

\section{Weighted bad-children estimates and survival}\label{sec:survival}
Proposition~\ref{prop:thin} controls thin inputs. To control a persistent
pair comparison, we first show that it forces two short vectors on an
interval of controlled length. The resulting occupation estimate can
then be integrated using the area-fiber bounds.

For a fixed height word $\mathbf W=(W_0,\ldots,W_n)$, write
$\hat Q_{j}:=\hat Q(W_j)=W_j^{\hat \gamma}$ for $0\le j\le n$, where
$\hat \gamma=(\gamma-1)/2$.
Similarly, set $R_j:=W_j^{\eta \hat \gamma}$,
$\sigma_j:=R_j^{-b_0}$, and $D_j:=R_j^{A+1}$.
These quantities are evaluated at the representative heights; their values at the actual slit lengths differ only by fixed factors.
\begin{lemma}\label{lem:occupation}
Suppose $\hat{Q}\ge1$ and $\lambda_1(z_i)\ge \epsilon$, where
$0<\epsilon\le1$, and the pair $(i,j)$ satisfies
\eqref{eq:persistent-comparison} for every $b\in[-1,1]$.
There is an interval of length $c/\hat{Q}$ in $[-1,1]$ on which
\begin{equation}\label{eq:occupation-thresholds}
\begin{aligned}
\lambda_1\bigl((a_{\frac12\log \hat{Q}}u_bz)_i\bigr)
&\le\epsilon_i:=\frac{C}{\epsilon\sqrt{\hat{Q}}},\\
\lambda_1\bigl((a_{\frac12\log \hat{Q}}u_bz)_j\bigr)
&\le\epsilon_j:=\frac{CD^{5/2}}{\epsilon\sqrt{\hat{Q}}}.
\end{aligned}
\end{equation}
\end{lemma}

\begin{proof}
Standard lattice reduction bounds the covering radius of $z_i$ by
$C_0/\epsilon$. Choose $(x,y)\in z_i$ within $C_0/\epsilon$ of
$(-B/(2\epsilon),B/\epsilon)$, where $B>16C_0$ is fixed.
Then $y\asymp\epsilon^{-1}$ and $b_*=-x/y\in[1/4,3/4]$.
For $|b-b_*|\le c/\hat{Q}$,
\begin{equation*}
a_{\frac12\log \hat{Q}}u_b(x,y)
=\left(\sqrt{\hat{Q}}\,y(b-b_*),\frac{y}{\sqrt{\hat{Q}}}\right),
\end{equation*}
whose length is at most $C/(\epsilon\sqrt{\hat{Q}})$.
Take $c\le1/8$ so that the interval lies in $[-1,1]$.
This nonzero witness proves the first bound; the persistent comparison
proves the second. Adjusting $c$ gives the stated interval length.
\end{proof}

\begin{proposition}\label{prop:bad}
For a height word $\mathbf W$ of depth $n\ge1$,
\begin{equation}\label{eq:weighted-bad}
\nu_{n,\mathbf W}(\mathcal B)
\le C\delta_{n-1}^g E_n(1+W_0^{r\varepsilon})
\left(\sigma_n^2+
\frac{D_n^5\sigma_n^{-4}}{\hat Q_{n}}\right).
\end{equation}
Here $E_n$ is defined in Lemma~\ref{lem:history}; all thresholds are
evaluated at the representative height $W_n$, with fixed comparison
constants covering the actual heights.
\end{proposition}

\begin{proof}
Proposition~\ref{prop:thin} bounds inputs with
$\min_k\lambda_1(\Lambda_k)<\sigma_n$ by the first term.
For the remaining inputs, fix a pair $i<j$ satisfying the persistent
comparison. The bounded shear and rotation defining $z$ give
$\lambda_1(z_i)\ge c\sigma_n$.
Write $\ell_k(b)=\lambda_1((a_{\frac12\log \hat{Q}}u_bz)_k)$.
Lemma~\ref{lem:occupation} gives
\begin{equation*}
\mathbf1_{\mathrm{comparison}}
\le C \hat Q_{n}\int_{-1}^1
\mathbf1_{\{\ell_i(b)\le\epsilon_i,\,
\ell_j(b)\le\epsilon_j\}}\,db,
\end{equation*}
where the thresholds in \eqref{eq:occupation-thresholds}, with
$\epsilon=c\sigma_n$, can be chosen uniformly on the height word.
At each fixed time, Lemma~\ref{lem:two-area} gives the area-volume
saving $C\epsilon_i^2\epsilon_j^2$.
The shift $K-1/Q$ depends only on the fixed slit output and is allowed
in that estimate.

We transfer this occupation estimate to the radial cone used in
Lemma~\ref{lem:fiber}. For a dilation $Z\mapsto\rho Z$ in its fixed
positive radial interval, retain the same integer history and put
\[
\kappa_\rho=\rho^{(1+\gamma)/2}
\sqrt{\frac{\log(\rho W)}{\log W}}.
\]
Area normalization cancels the scalar $\rho$, leaving
$\Lambda(\rho Z)=a_{-\log\kappa_\rho}\Lambda(Z)$.
The identity $k a_{-\log\kappa_\rho}=a_{\log\kappa_\rho}k$
gives the change of time
\begin{equation}\label{eq:radial-time}
b_\rho=\kappa_\rho^2(b+s)-s_\rho,
\qquad s=K-1/Q(W),\qquad s_\rho=K-1/Q(\rho W).
\end{equation}
Indeed, $u_c a_v=a_vu_{ce^{-2v}}$. If
$F_b(Z)=a_{\frac12\log \hat{Q}(W)}u_bz(Z)$, then
\begin{equation}\label{eq:radial-orbit}
F_{b_\rho}(\rho Z)
=a_{\frac12\log(\hat{Q}(\rho W)/\hat{Q}(W))+\log\kappa_\rho}
F_b(Z).
\end{equation}
The diagonal factor and its inverse have uniformly bounded norms,
since $\hat{Q}(\rho W)/\hat{Q}(W)=\rho^{\hat{\gamma}}$.
Also $db_\rho=\kappa_\rho^2\,db$, with a uniformly bounded positive
Jacobian. Thus the transformed time intervals lie in one fixed bounded
interval $J_*$, and the original short-vector event is contained in
the transformed event after a fixed enlargement of both thresholds.

The full period space has real dimension $4g+2$, so radial
disintegration contributes the bounded positive density
$\rho^{4g+1}\,d\rho$ on this radial interval.
Apply \eqref{eq:radial-time} separately at each cone point before
Tonelli's theorem. This remains valid when the radial parameter depends
on the area coordinates of the subsequent fiber description.
Dilation changes the geometric windows only by fixed factors; the
containing fiber domains absorb these changes. It need not preserve
the normality tests or candidate sets. For each history we use the
supremum of its incoming weight over the original area-one domain,
as in Lemma~\ref{lem:fiber}.

At each fixed transformed time, the full-period volume estimate in
Lemma~\ref{lem:fiber}, combined with Lemma~\ref{lem:two-area}, is
\[
C\delta_{n-1}^gE_n(1+W_0^{r\varepsilon})
\epsilon_i^2\epsilon_j^2,
\]
uniformly on $J_*$. The normalization of the root measure and the
sector summation are included in that estimate. Integrating over
$J_*$ and using
\begin{equation*}
\hat Q_{n}\epsilon_i^2\epsilon_j^2
\le C \hat Q_{n}
\frac{\sigma_n^{-2}}{\hat Q_{n}}
\frac{D_n^5\sigma_n^{-2}}{\hat Q_{n}}
=\frac{CD_n^5\sigma_n^{-4}}{\hat Q_{n}}
\end{equation*}
bounds the arrival mass for this pair by the second term of
\eqref{eq:weighted-bad}. Summing over the $\binom g2$ pairs proves
the proposition with the original history weights.
\end{proof}

\begin{proof}[Proof of Theorem~\ref{thm:quantitative}]
Fix the constants $A,\kappa,P,d_1$ in the preceding estimates.
After fixing $g$ and a sufficiently small $t>0$, choose $b_0$, then
$\eta$, and finally $\varepsilon$ so that
\begin{equation}\label{eq:parameters}
\begin{gathered}
0<b_0<\frac1{4A},\qquad C_p=5(A+1),\\
0<\eta<\min\left\{\frac{\log2}{4A},
\frac1{4(C_p+4b_0)}\right\},\\
m_*:=\min\{2b_0\eta,\,1-\eta(C_p+4b_0)\}>0,\\
0<r\varepsilon<\min\left\{\frac t2,
\frac{tm_*}{16},\frac{\kappa\eta t}{8(d_1+1)}\right\}.
\end{gathered}
\end{equation}
These choices are independent of the height and permit arbitrarily
small positive $\varepsilon$. Recall that $\hat\gamma=(\gamma-1)/2\ge t/2$,
$R=W^{\eta \hat{\gamma}}$, $\sigma=R^{-b_0}$, and $D=R^{A+1}$.
Thus
\[
R_{\log}=\frac{\eta \hat{\gamma}\log W}{\log2},\qquad
T=q\log W+\tfrac12\log\delta^{-1}.
\]
For large $W$, the bound on $\eta$ gives
$A R_{\log}<q\log W/4$, leaving a uniform margin in the time
condition of Theorem~\ref{thm:effective}, including its integer-time
formulation. Moreover,
$\inj(z)^{-1}\le C\sigma^{-2}=C R^{2b_0}=o(R^{1/A})$ at normal
inputs. Lemma~\ref{lem:obstruction} therefore applies to each such input.

To compare the effective error with the shrinking mean, observe that
\begin{equation*}
\delta^{-(P+g)}\tau^{-(d_1+1)}R^{-\kappa}
=O\!\left((\log W)^{P+g}
W^{(d_1+1)r\varepsilon-\kappa\eta \hat{\gamma}}\right)=o(1).
\end{equation*}
Indeed, $(d_1+1)r\varepsilon<\kappa\eta t/8\le\kappa\eta \hat{\gamma}/4$.
Proposition~\ref{prop:count} now gives the branching bound uniformly
at every normal input above a fixed height threshold. Increasing $W_0$
puts all subsequent inputs above that threshold.

The two savings in \eqref{eq:weighted-bad} are
\begin{equation*}
\sigma_n^2=W_n^{-2b_0\eta \hat{\gamma}},\qquad
\frac{D_n^5\sigma_n^{-4}}{\hat Q_{n}}
=W_n^{-q[1-\eta(C_p+4b_0)]}.
\end{equation*}
Put $\beta=\hat\gamma m_*>0$. Both are at most $W_n^{-\beta}$.
Summing over the at most $C^n$ height words, absorbing this factor
into $E_n$, and using $\delta_{n-1}\le1$, we obtain
\begin{equation}\label{eq:depth-loss}
L_n\le E_n(1+W_0^{r\varepsilon})e^{-\beta h_n},
\qquad
E_n=\exp(Cn^2+Cn\log\log W_0+Cn),
\end{equation}
where $h_n$ is the least logarithmic height at depth $n$.
Write $H_0=\log W_0$. The recurrence
$|h_{j+1}-\gamma h_j|\le C_h$ along every history gives
\begin{equation}\label{eq:height-growth}
h_n\ge\gamma^n H_0-C_h\frac{\gamma^n-1}{\gamma-1}
\ge\frac34\gamma^nH_0
\end{equation}
once $H_0\ge4C_h/(\gamma-1)$.
For fixed $\gamma>1$,
$\sup_{n\ge1}(n^2+n)/\gamma^n<\infty$ and
$\log H_0/H_0\to0$. Hence we may further increase $H_0$ so that
\[
\log E_n+\log2\le\frac\beta8\gamma^nH_0
\qquad(n\ge1).
\]
Since $r\varepsilon<tm_*/16\le\beta/8$, it follows uniformly in
$n\ge1$ that
\[
\log E_n+\log(1+W_0^{r\varepsilon})
\le\frac\beta4\gamma^nH_0.
\]
Combining this with \eqref{eq:depth-loss}--\eqref{eq:height-growth}
yields
\begin{equation}\label{eq:superexponential-loss}
L_n\le\exp(-\tfrac\beta2\gamma^nH_0),\qquad n\ge1.
\end{equation}
This proves the asserted positive-depth bound with $c_* =\beta/2$.

At the root, the normalized lattices lie in fixed compact charts, so
thin inputs are absent for large $W_0$. Conditional on the slit and
component areas, the root measure has uniformly bounded density
relative to product Haar measure. This bound persists under the common
shear and rotation defining $z$. Apply Lemma~\ref{lem:occupation}
with a fixed $\epsilon>0$. At each fixed time, Haar invariance and the
primitive Siegel formula bound the two-short-vector event by
$C\epsilon_i^2\epsilon_j^2$. Integrating the occupation bound and
summing over pairs gives
\begin{equation}\label{eq:root-loss}
L_0\le C_g\frac{D_0^5}{\hat Q_{0}}
=C_g W_0^{-q[1-5\eta(A+1)]}\longrightarrow0.
\end{equation}
The exponent is positive by \eqref{eq:parameters}; dyadic variation
of the root height changes only the fixed constant.

The bound \eqref{eq:superexponential-loss} is summable in $n$, and
its sum tends to zero as $W_0\to\infty$. Together with
\eqref{eq:root-loss}, this proves the asserted total-loss estimate.
Finally, \eqref{eq:conservation} gives
\[
\int_\Omega s\,d\mu_\Omega=1-\sum_{n\ge0}L_n>0
\]
for sufficiently large $W_0$. Since $0\le s\le1$, the roots with
$s(Z)>0$ form a set of positive measure. The area and angle budgets
established above hold along every retained history, completing the
proof.
\end{proof}

\section{From surviving paths to almost every surface}\label{sec:directions}

We first fix a root and construct a measure on its surviving paths. The
resulting pointwise dimension estimate will then be promoted to an
almost-everywhere statement in the ambient connected component.

\begin{lemma}[The surviving path measure]\label{lem:path-measure}
For a root $Z$, let $s_N(Z)$ be the total original mass of paths that
remain normal through depth $N$. The limit
$s(Z)=\lim_{N\to\infty}s_N(Z)$ is Borel measurable in $Z$. There is a
finite Borel measure $\beta_Z$ on the infinite surviving paths, of total
mass $s(Z)$, such that
\begin{equation}\label{eq:path-cylinder-mass}
  \beta_Z([S])\le p(S)=\prod_{j<|S|}b(S_j)^{-1}
\end{equation}
for every normal vertex $S$, where $[S]$ is its path cylinder and $S_j$
are its strict ancestors. For sufficiently large $W_0$, the set of roots
with $s(Z)>0$ has positive period measure.
\end{lemma}

\begin{proof}
For a nonnormal root, set $s_N(Z)=s(Z)=0$ and $\beta_Z=0$. Otherwise, for
a normal vertex $S$ and $N\ge |S|$, let $q_N(S)$ be the sum of the original
incoming weights of its normal descendants at depth $N$. The sequence
$q_N(S)$ decreases and is bounded by $p(S)$. Its limit $q(S)$ satisfies
\[
  q(S)=\sum_{S'\text{ child of }S}q(S'),
\]
with failed children assigned mass zero, because every vertex has
finitely many children. These consistent cylinder masses define a finite
Borel measure on the boundary of the surviving tree. Its total mass is
$q(Z)=s(Z)$, and $q(S)\le p(S)$ gives
\eqref{eq:path-cylinder-mass}.

The finite-depth masses are Borel by the construction in
Section~\ref{sec:fibers}, so their decreasing limit $s$ is Borel.
Theorem~\ref{thm:quantitative} and \eqref{eq:conservation} give
$\int_\Omega s\,d\mu_\Omega>0$ for sufficiently large $W_0$. Since
$0\le s\le1$, the set $\{s>0\}$ has positive measure.
\end{proof}

\begin{lemma}[Invariant regions from a surviving path]
\label{lem:invariant-region}
Every infinite surviving path has a limiting direction $\theta$. On its
fixed surface, the flow in direction $\theta$ admits an invariant
measurable set $F$ with $0<\area(F)<1$.
\end{lemma}

\begin{proof}
Let $w_n$ be the slit, $W_n=|w_n|$, and $F_n=T_1^n$ the first torus
region at depth $n$. The move estimates and the height recurrence give
\begin{equation}\label{eq:surviving-path-estimates}
\begin{gathered}
  \sum_n\area(F_n\triangle F_{n+1})<\infty,
  \qquad |\arg w_n-\theta|\le C\delta_nW_n^{-1-\gamma},\\
  a\le\area(F_n)\le1-(g-1)a.
\end{gathered}
\end{equation}
Here the existence of $\theta$ and the angular estimate follow by summing
the turning bounds, using $W_{n+1}\asymp W_n^\gamma$ and $\gamma>1$.
Summability of the symmetric differences implies that
$\mathbf1_{F_n}$ converges almost everywhere and in $L^1$ to the
indicator of a measurable set $F$. The area bounds pass to the limit.

The boundary of $F_n$ consists of two slit banks. The area swept out by
these banks under a flow segment of length $|u|$ is at most
$2|u|\,|w_n\times e_\theta|$, where $e_\theta$ is a unit vector in
direction $\theta$. Thus
\begin{equation*}
  \area(F_n\triangle\phi_\theta^uF_n)
  \le2|u|\,|w_n\times e_\theta|
  \le C|u|\delta_nW_n^{-\gamma}\longrightarrow0.
\end{equation*}
Trajectories meeting a zero form an area-null set. Since the flow
preserves area, for every fixed $u$ we obtain
\[
  \area(F\triangle\phi_\theta^uF)
  \le2\area(F\triangle F_n) +\area(F_n\triangle\phi_\theta^uF_n)\longrightarrow0.
\]
Hence $F$ is invariant modulo null sets and has area strictly between
zero and one.
\end{proof}

\begin{proposition}
\label{prop:dimension}
Fix admissible $t,\varepsilon$ and a root $Z$ with $s(Z)>0$. On its surface,
\begin{equation}\label{eq:pointwise-dimension}
  \dim_H\NE(X,\omega)
  \ge\frac{t}{(\gamma-1)(2\gamma-t)},
  \qquad \gamma=(1+t)(1+\varepsilon).
\end{equation}
\end{proposition}

\begin{proof}
Push $\beta_Z$ forward by the limiting-direction map of
Lemma~\ref{lem:invariant-region}, and denote the resulting measure by
$\nu$. It has mass $s(Z)>0$ and is supported on $\NE(X,\omega)$. The
direction map is continuous: agreement through depth $n$ confines the
limiting directions to an interval whose diameter tends uniformly to
zero. In particular, its image is compact.

Choose $a_-,a_+>0$ such that
\[
  a_-\gamma^n\le\log W_n\le a_+\gamma^n
\]
throughout the tree. The recurrence
$\log W_{n+1}=\gamma\log W_n+O(1)$ allows
$a_\pm=\log W_0+O((\gamma-1)^{-1})$. Take the root tall enough that
$(1+\gamma)a_->2a_+$. By \eqref{eq:surviving-path-estimates}, an interval
centered at $[w_n]$ of radius $cW_n^{-1-\gamma}$ contains all descendant
directions. These intervals are nested for sufficiently large $W_0$:
the turn to a child is smaller than the parent radius by a factor
$O(\delta_n)$, and the child's radius is smaller still. The sibling
separation in Lemma~\ref{lem:children} and the inequality
$(1+\gamma)a_->2a_+$ give a gap
\[
  \zeta_n=c\exp(-2a_+\gamma^n)
\]
between distinct depth-$n$ intervals. Branches that separate earlier
have larger gaps.

Every normal vertex at depth $n$ has at least $m_n$ candidates, where
\begin{equation*}
  m_n=\left\lfloor
      \frac{c\exp(ta_-\gamma^n)}{(a_+\gamma^n)^g}
      \right\rfloor\ge2,
  \qquad M_n=\prod_{j<n}m_j.
\end{equation*}
Lemma~\ref{lem:path-measure} bounds the mass of each surviving depth-$n$
cylinder by $M_n^{-1}$. For $\zeta_{n+1}\le z<\zeta_n$, a ball of
radius $z$ meets at most a fixed number of depth-$n$ intervals and at
most $Cz/\zeta_{n+1}$ depth-$(n+1)$ intervals. Consequently, for
$0<s<1$,
\begin{equation}\label{eq:frostman-interpolation}
  \nu(B(\theta,z))
  \le C\min\left\{M_n^{-1},
                   \frac{z}{M_{n+1}\zeta_{n+1}}\right\}
  \le\frac{Cz^s}{M_n(m_n\zeta_{n+1})^s}.
\end{equation}
The last inequality uses $M_{n+1}=M_nm_n$ and
$\min\{A,B\}\le A^{1-s}B^s$. For this fixed root and fixed parameters,
$\log m_n=ta_-\gamma^n+O(n+1)$, and hence
\begin{equation*}
  \log\bigl[M_n(m_n\zeta_{n+1})^s\bigr]
  \ge\left(\frac{ta_-}{\gamma-1}
           -s(2\gamma a_+-ta_-)\right)\gamma^n-O(n^2).
\end{equation*}
If $ s<\frac{ta_-}{(\gamma-1)(2\gamma a_+-ta_-)}$,
the right-hand side tends to infinity. After enlarging the constant to
cover finitely many scales, \eqref{eq:frostman-interpolation} gives a
Frostman bound for $\nu/s(Z)$. The mass distribution principle therefore
gives the corresponding dimension lower bound for its compact support.

At every depth, some cylinder has positive surviving mass. Restricting
to such a cylinder $[S]$ and dividing by its positive original incoming
mass $p(S)$ produces a subtree with unit incoming mass and positive
surviving mass $\beta_Z([S])/p(S)$. The same cylinder bounds apply in this
subtree. As its root height tends to infinity, the corresponding ratio
$a_-/a_+$ tends to one while $t,\varepsilon$ remain fixed. Taking the
supremum of these lower bounds proves \eqref{eq:pointwise-dimension}.
\end{proof}

We now use the standard finiteness and ergodicity
\cite{Masur1982,Veech1982,Veech1986}. Recall also that
$\NE(Y)\subseteq\NUE(Y)$ and that Masur's bound gives
$\dim_H\NUE(Y)\le1/2$ for every compact translation surface of genus at
least two \cite{Masur1992}.

\begin{proposition}
\label{prop:promotion}
Let $\mathcal C$ be a connected component of an area-one stratum of
genus at least two, with normalized measure $\mu$. Suppose
that $s_j\to1/2$ and, for every $j$, there is a $\mu$-measurable set
$E_j\subset\mathcal C$ of positive measure such that
$\dim_H\NE(Y)\ge s_j$ for every $Y\in E_j$. Then
\begin{equation}\label{eq:almost-everywhere-equality}
  \dim_H\NE(Y)=\dim_H\NUE(Y)=\frac12
\end{equation}
for $\mu$-almost every $Y\in\mathcal C$.
\end{proposition}

\begin{proof}
For $A\in\SL_2(\R)$, the induced affine map sends directional flows to
directional flows, with a constant time change in each direction, and
preserves flat area. Thus
\[
  \NE(AY)=D_A(\NE(Y)),\qquad D_A([v])=[Av].
\]
The projective map $D_A$ is a smooth circle diffeomorphism and hence
bi-Lipschitz. It follows that $\dim_H\NE$ is invariant under
$\SL_2(\R)$.

By inner regularity, choose compact sets $K_j\subset E_j$ with
$\mu(K_j)>0$. Each saturation
$\SL_2(\R)K_j$ is a countable union of compact sets, because
$\SL_2(\R)$ is $\sigma$-compact and its action is continuous. It is
therefore Borel, invariant, and of positive measure, so ergodicity gives
$\mu(\SL_2(\R)K_j)=1$. On this saturation, $\dim_H\NE\ge s_j$.
The countable intersection of these saturations has full measure and
satisfies $\dim_H\NE\ge1/2$. The inclusion
$\NE\subseteq\NUE$ and Masur's upper bound yield
\eqref{eq:almost-everywhere-equality}.
\end{proof}

\begin{proof}[Proof of Theorem~\ref{thm:main}]
Fix $g\ge3$. Choose a sequence $t_j\downarrow0$ and, for each $j$,
choose $0<\varepsilon_j\le t_j^3$ satisfying \eqref{eq:parameters}.
Theorem~\ref{thm:quantitative} and Lemma~\ref{lem:path-measure} give a
positive-measure Borel set of roots with positive surviving mass for
each parameter pair. By Lemma~\ref{lem:component}, the cyclic period
maps are locally nonsingular in $\mathcal C_g$. A countable cover by
injective local charts therefore contains a chart that maps a
positive-measure part of this root set to a positive 
measure set $E_j\subset\mathcal C_g$.

Put $\gamma_j=(1+t_j)(1+\varepsilon_j)$. Proposition~\ref{prop:dimension}
gives, for every $Y\in E_j$,
\[
  \dim_H\NE(Y)\ge s_j
  :=\frac{t_j}{(\gamma_j-1)(2\gamma_j-t_j)}
  \longrightarrow\frac12,
\]
since $\gamma_j-1=t_j+O(t_j^3)$ and $2\gamma_j-t_j=2+O(t_j)$.
Proposition~\ref{prop:promotion}, applied in the fixed component
$\mathcal C_g$, proves the theorem.
\end{proof}

\begin{corollary}\label{cor:components}
If $g=2$, $g=3$, or $g\ge4$ is even, then for almost every
area-one surface in every connected component of
$\mathcal H(g-1,g-1)$, the set of non-uniquely ergodic directions has
Hausdorff dimension $1/2$. For odd $g\ge5$, the same conclusion holds
in the odd-spin and hyperelliptic components.
\end{corollary}

\begin{proof}
Theorem~\ref{thm:main} treats $\mathcal C_g$ for $g\ge3$, while
Chaika--Masur \cite[Theorem~1.5]{ChaikaMasur} gives the conclusion in
the hyperelliptic component for every $g\ge2$. Theorem~\ref{thm:KZ-specialized} gives the
stated cases.
\end{proof}

\end{document}